\documentclass[aos,MSNbibl,nameyear,dvips]{arximspdf}

\doi{10.1214/12-AOS1013} %kopijuoti is PTS
\volume{40}
\issue{3}
\pubyear{2012}
\firstpage{1550}
\lastpage{1577}

\makeatletter
\newcommand{\eqref}[1]{(\ref{#1})}
\newcommand{\N}{\mathcal{N}}
\newcommand{\St}{\mathcal{S}t}
\newproclaim{example}{Example}%[section]
\newtheorem{crit}{Criterion}
\newproclaim{deft}{Definition}
\newproclaim{argg}{Argument}
\newtheorem{lem}{Lemma}
\newtheorem{prop}{Proposition}
\newtheorem{res}{Result}
\newtheorem{cor}{Corollary}
\makeatother

\begin{document}
\begin{frontmatter}

\title{Criteria for Bayesian model choice with application to variable selection\thanksref{T1}}
\runtitle{Criteria for Bayesian model choice}

\thankstext{T1}{Supported in part by the Spanish Ministry of Education and Science
Grant MTM2010-19528, and by NSF Grants DMS-06-35449, DMS-07-57549-001 and DMS-10-07773.}

\begin{aug}
\author[A]{\fnms{M.~J.} \snm{Bayarri}\ead[label=e1]{susie.bayarri@uv.es}},
\author[B]{\fnms{J.~O.} \snm{Berger}\ead[label=e2]{berger@stat.duke.edu}},
\author[C]{\fnms{A.} \snm{Forte}\corref{}\ead[label=e3]{forte@uji.es}}
\and
\author[D]{\fnms{G.} \snm{Garc\'ia-Donato}\ead[label=e4]{Gonzalo.GarciaDonato@uclm.es}}
\runauthor{M. J. Bayarri et al.}
\affiliation{Universitat de Val\`encia, Duke University, Universitat
Jaume I and Universidad de Castilla-La Mancha}
\address[A]{M.~J. Bayarri\\
Department of Mathematics\\
Universitat de Val\`encia\\
Valencia, Valencia\\
Spain\\
\printead{e1}} %adresu isvedimo komanda gale!
\address[B]{J.~O. Berger\\
Department of Statistics\\
Duke University\\
Durham, North Carolina\\
USA\\
\printead{e2}\hspace*{44pt}}
\address[C]{A. Forte\\
Department of Economics\\
Universitat Jaume I\\
Castell\'{o}n, Valencia\\
Spain\\
\printead{e3}}
\address[D]{G. Garc\'ia-Donato\\
Department of
Economic Analysis\\
\quad and Finance\\
Universidad de Castilla-La Mancha\\
Albacete, Castilla-La Mancha\\
Spain\\
\printead{e4}}
\end{aug}

% HISTORY:
\received{\smonth{7} \syear{2011}}
\revised{\smonth{4} \syear{2012}}

% ABSTRACT
%
\begin{abstract}
In objective Bayesian model selection, no single criterion has emerged as
dominant in defining objective prior distributions. Indeed, many criteria
have been separately proposed and utilized to propose differing prior choices.
We first formalize the most general and compelling of the various criteria
that have been suggested, together with a new criterion. We then
illustrate the potential of these criteria in determining objective
model selection priors by considering their application to the
problem of variable selection in normal linear models. This results in a
new model selection objective prior with a
number of compelling properties.
\end{abstract}

% KEYWORDS
%


\end{frontmatter}

%-----------------------------Introduction
%s1 #&#
\section{Introduction}\label{sec1}
%--------------------------------------------
%s1.1 #&#
\subsection{Background}
\label{secprelim}

A key feature of Bayesian model selection, when the models have
differing dimensions and noncommon parameters, is that results are typically
highly sensitive to the choice of priors for the noncommon parameters,
and, unlike the scenario for
estimation, this sensitivity does not vanish as the sample size grows; see
\citet{KassRaf95}, \citet{BerPer01}. Furthermore,
improper priors cannot typically be used for noncommon parameters,
nor can ``vague proper priors'' (see the above references, e.g., and the
brief discussion in Section~\ref{secbasic}), ruling out use of the
main tools developed in objective Bayesian estimation theory.

Because of the difficulty in assessing subjective priors for numerous
models, there have been many efforts (over more than 30 years) to
develop ``conventional'' or
``objective'' priors for model selection; we will term these ``objective
model selection priors,''
the word objective simply meant to indicate that they are not
subjective priors, and
are chosen conventionally based on the models being considered. A few
of the many references most related to this paper are
\citet{Jef61},
\citeauthor{ZellSiow80} (\citeyear{ZellSiow80,ZellSiow84}),
\citet{LauIbra95},
\citet{KassWass95},
\citet{BerPer96},
\citet{MorBerRac98},
\citet{DeSpe99},
\citet{PerBer02},
\citet{BayGar08},
\citet{liang08},
\citet{CuiGeo08},
\citet{MarGeo08},
\citet{MarStr10}.

For the most part, these efforts were started with a good idea which
was used
to develop the priors, and then the behavior of the priors was studied.
Yet, in spite of the apparent success of many of these methods, there has
been no agreement as to which are most appealing or most successful.

This lack of progress in reaching a consensus on objective priors for
model selection resulted in our approaching the problem from a different
direction, namely, formally formulating the various criteria that have
been deemed essential for model selection priors (such as consistency
of the resulting procedure), and seeing if these
criteria can essentially determine the priors.

The criteria are stated for general model selection problems in
Section~\ref{secCriteria},
which also discusses their historical antecedents. To illustrate that
application of
the criteria can largely determine model selection priors, we turn to a
specific problem
in Section~\ref{priors}---variable
selection in normal linear models. The resulting priors for variable
selection are new
and result in closed form Bayes factors; for those primarily interested
in the methodology itself,
the resulting priors and Bayes factors are given in Section~\ref{secsummary}.

%s1.2 #&#
\subsection{Notation}
Let $\mathbf{y}$ be a data vector of size $n$ from one of the models
%
%e1 #&#
\begin{equation}
\label{eqTheGenProb} M_0\dvtx f_0(\mathbf{y}\mid\bolds{
\alpha}),\qquad M_i\dvtx f_i(\mathbf{y}\mid\bolds{\alpha},
\bolds{\beta }_i),\qquad i=1,2,\ldots, N-1 ,
\end{equation}
where $\bolds{\alpha}$ and the $\bolds{\beta}_i$ are
unknown model parameters, the
latter having dimension~$k_i$. $M_0$ will be called the null model and
is nested in all of the considered models.

Under the null model, the prior is $\pi_0(\bolds{\alpha})$;
under model
$M_i$, and without loss of generality, we express the model selection
prior as
\[
\pi_i(\bolds{\alpha},\bolds{\beta}_i)=
\pi_i(\bolds {\alpha}) \pi_i(\bolds{
\beta}_i\mid\bolds{\alpha}).
\]
Note that the parameter $\bolds{\alpha}$ occurs in all of the
models, so that
$\bolds{\alpha}$ is usually referred to as the \textit{common}
parameter; the $\bolds{\beta}_i$ are called \textit{model specific}
parameters.
%see e.g.][] {Lindley00,BerPer01}

%%%-------------------------------------------------------------------------------------

Assuming that one of the entertained models is true, the posterior
probability of each of the models $M_i$ can be written in the
convenient form
%
%e2 #&#
\begin{equation}
\label{postprob} \operatorname{Pr}(M_i\mid\mathbf{y})=\frac{B_{i0}}{1+ (\sum_{j=1}^{{N}-1}
B_{{j}{0}}
P_{{j}{0}}  )} ,
\end{equation}
where $P_{{j}{0}}$ is the prior odds $P_{{j}{0}}=\operatorname{Pr}(M_j)/\operatorname{Pr}(M_0)$, with
$\operatorname{Pr}(M_j)$ being the prior probability of model $M_j$, and $B_{{j}{0}}$
is the Bayes factor of model $M_j$ to~$M_0$ defined by
%
%e3 #&#
\begin{equation}
\label{eqmarginal} B_{{j}{0}}= \frac{m_j(\mathbf{y})}{m_0(\mathbf{y})}\qquad \mbox{with }
m_j(\mathbf{y})=\int f_{j}(\mathbf{y}\mid\bolds{\alpha },
\bolds{\beta}_i) %
{\pi_{j}(\bolds{\alpha},
\bolds{\beta}_{j})} %
\,d\bolds{\alpha} \,d\bolds{
\beta}_{j}
\end{equation}
and $m_0(\mathbf{y})=\int f_{0}(\mathbf{y}\mid\bolds
{\alpha}){\pi_{0}(\bolds{\alpha})}
\,d\bolds{\alpha}$ being the marginal likelihoods of model~$M_j$
and $M_0$
corresponding to the model prior densities
$\pi_{j}(\bolds{\alpha},  \bolds{\beta}_{j})$ and $\pi_{0}(\bolds{\alpha})$. [Any
model could serve as the base model for computation of the Bayes
factors in (\ref{postprob}), but use of the null model is common and
convenient.] The focus in this paper is on choice of model priors $\pi_{0}(\bolds{\alpha})$ and $\pi_{j}(\bolds{\alpha},
\bolds{\beta}_{j})$.

%Choice of the model prior probabilities $\operatorname{Pr}(M_i)$ is also
%important, with equal prior probabilities often not being
%suitable
%but this choice is not germane to this paper.

%%%------------------------------------------------------------------------------------

%-------------------------------Criteria for objective model selection
%priors
%s2 #&#
\section{Criteria for objective model selection priors}\label{secCriteria}

%---------------------------------------------------------------------
%s2.1 #&#
\subsection{Introduction}
\label{secfcrit}

The arguments concerning prior choice in testing and model selection in
\citet{Jef61} are often called Jeffreys's \textit{desiderata} [see \citet
{BerPer01}] and are the precursors to the
criteria developed herein. [\citet{Robetal09}, is a comprehensive and
modern review of Jeffreys's book.]
These and related ideas have been repeatedly used to evaluate or guide
development of objective model priors;
see, for example, \citet{BerPer01}, \citet{BayGar08}, \citet{liang08} and \citet
{Forte11}. We group the criteria into four classes: basic, consistency
criteria, predictive matching criteria and invariance criteria.

%----------------------------
%s2.2 #&#
\subsection{Basic criteria}
\label{secbasic}

As mentioned in the \hyperref[sec1]{Introduction} priors for the noncommon parameters
$\bolds{\beta}_i$ should be proper, because they only occur in
the numerator
of the Bayes factors $B_{i0}$, and hence, if using an improper prior,
the arbitrary constant for the improper prior would not cancel, making
$B_{i0}$ ill defined. There have been various efforts to use improper
priors and define a meaningful scaling [\citet{GhosSamanta2002}, \citet{SpSm82}];
and other methods have been proposed that can be interpreted as
implicitly scaling the improper prior Bayes factor [see details and
references in \citet{BayGar08}], but we are restricting consideration
here to real Bayesian procedures.

Similarly, vague proper priors cannot be used in determining the
$B_{i0}$, since the arbitrary scale of vagueness appears as a
multiplicative term in the Bayes factor, again rendering the Bayes
factor arbitrary. Thus we have:

\begin{crit}[(Basic)]
Each conditional prior $\pi_i(\bolds{\beta}_i\mid
\bolds{\alpha})$ must be
proper (integrating to one) and cannot be arbitrarily vague in the
sense of almost all of its mass being outside any believable compact set.
\end{crit}

%s2.3 #&#
\subsection{Consistency criteria}
\label{secconsiscrit}
Following \citet{liang08}, we consider two primary consistency criteria---model
selection consistency and information
consistency:\vadjust{\goodbreak}

\begin{crit}[(Model selection consistency)]
If data $\mathbf{y}$ have been generated by $M_i$, then the posterior
probability of $M_i$ should converge to 1 as the sample size $n
\rightarrow\infty$.\vspace*{-1pt}
\end{crit}

Model selection consistency is not particularly controversial, although
it can be argued that the true model is never one of the entertained
models, so that the criterion is vacuous. Still, it would be
philosophically troubling to be in a situation with infinite data
generated from one of the models being considered, and not choosing the
correct model. A number of recent references concerning this criterion
are \citet{FLS01}, \citet{BerGhoMuk03}, \citet{liang08},
\citet{CasGirMarMor09}, \citet{GuoSpe09}.\vspace*{-1pt}

\begin{crit}[(Information consistency)]
For any model $M_i$, if $\{ \mathbf{y}_m, m=1, \ldots\}$ is a
sequence of
data vectors of fixed size such that, as $m \rightarrow
\infty$,
%
%e4 #&#
\begin{equation}
\label{eqlr} \Lambda_{i0}(\mathbf{y}_m)=\frac{\sup_{\bolds{\alpha
},\bolds{\beta}_{i}}f_{i}(\mathbf{y}_m\mid
\bolds{\alpha}, \bolds{\beta}_{i})}{\sup_{\bolds
{\alpha}}f_{0}(\mathbf{y}_m\mid\bolds{\alpha})}
\rightarrow\infty  \qquad\mbox{then } B_{i0}(\mathbf {y}_m)
\rightarrow \infty .\vspace*{-1pt}
\end{equation}
\end{crit}

In normal linear models, this is equivalent to saying that, if one
considers a sequence of data vectors for which the corresponding $F$
(or $t$) statistic goes to infinity, then the Bayes factor should also
do so for this sequence. \citet{Jef61} used this argument to justify a
Cauchy prior in testing that a normal mean is zero, and the argument
has also been highlighted in \citet{BerPer01}, \citet{BayGar08}, \citet{liang08}. One can
construct examples in which a real Bayesian answer violates information
consistency, but the examples are based on very small sample sizes and
priors with extremely flat tails. Furthermore, violation of information
consistency would place frequentists and Bayesians in a particularly
troubling conflict, which many would view as unattractive.

A third type of consistency has been proposed to address the fact that
objective model selection priors typically depend on specific features
of the model, such as the sample size or the particular covariates
being considered.\vspace*{-1pt}

\begin{crit}[(Intrinsic prior consistency)]
Let $\pi_i(\bolds{\beta}_i \mid\bolds{\alpha}, n)$
denote the prior for the
model specific parameters of model $M_i$ with sample size $n$. Then, as
$n \rightarrow\infty$ and under suitable conditions on the evolution
of the model with~$n$, $\pi_i(\bolds{\beta}_i \mid\bolds
{\alpha}, n)$ should
converge to a proper prior $\pi_i(\bolds{\beta}_i \mid
\bolds{\alpha})$.\vspace*{-1pt}
\end{crit}

The idea here is that, while features of the model and sample size (and
possibly even data) frequently affect model selection priors, such
features should disappear for large $n$. If there is such a limiting
prior, it is called an \textit{intrinsic prior}; see \citet{BerPer01} for
extensive discussion and previous references. (Note that some have used
the phrase ``intrinsic prior'' to refer to specific priors arising from
a specific model selection method, but we use the term here generically.)\vadjust{\goodbreak}

%s2.4 #&#
\subsection{Predictive matching criteria}
\label{secpmc}

The most crucial aspect of objective model selection priors is that
they be appropriately ``matched'' across models of different dimensions.
Having a prior scale factor ``wrong'' by a factor of 2 does not matter
much in one dimension, but in 50 dimensions that becomes an error of
$2^{50}$ in the Bayes factor. There have been many efforts to achieve
such matching in model selection, including \citet{SpSm82},
\citet{Suzuki83},
\citet{LauIbra95},
\citet{GhosSamanta2002}.

The standard approach to predictive matching is modeled after \citet
{Jef61}. For example, Jeffreys defined a ``minimal sample size'' for
which one would logically be unable to discriminate between two
hypotheses, and argued that the prior distributions should be chosen to
then yield equal marginal likelihoods for the two hypotheses. Here is
an illustration of this type of argument, from \citet{BerPerVar98}.

\begin{example*}
Suppose one is comparing two location-scale models
\[
M_1\dvtx y \sim\frac{1}{\sigma} p_1 \biggl(
\frac{y-\mu}{\sigma
} \biggr) \quad\mbox{and}\quad M_2\dvtx y \sim\frac{1}{\sigma}
p_2 \biggl(\frac
{y-\mu}{\sigma} \biggr) .
\]
Intuitively, two independent observations $(y_1,y_2)$ should not allow
for discrimination between these models, since two observations only
allow setting of the center and scale of the distribution; there are no
``degrees of freedom'' left for model discrimination.
Now consider the choice of prior (for both models) $\pi(\mu,\sigma) =
1/\sigma$. It is shown in \citet{BerPerVar98} that
\begin{eqnarray*}
&&\int\frac{1}{\sigma^2} p_1 \biggl(\frac{y_1-\mu}{\sigma} \biggr)
p_1 \biggl(\frac{y_2-\mu}{\sigma} \biggr) \pi(\mu,\sigma) \,d\mu \,d\sigma
\\
&&\qquad= \int\frac{1}{\sigma^2} p_2 \biggl(\frac{y_1-\mu}{\sigma} \biggr)
p_2 \biggl(\frac{y_2-\mu}{\sigma} \biggr) \pi(\mu,\sigma) \,d\mu \,d
\sigma= \frac{1}{2|y_1-y_2|}
\end{eqnarray*}
for any pair of observations $y_1 \neq y_2$, so that the models would
be said to be predictively matched for all minimal samples.
The Bayes factor between the models is then obviously 1, agreeing with
the earlier intuition that a~minimal sample should not allow for model
discrimination.
\end{example*}

This argument was formalized by \citet{BerPer01} as follows.

%de1 #&#
\begin{deft} The model/prior pairs $\{M_i,\pi_i\}$ and $\{M_j,\pi_j\}
$ are \textit{predictive matching at sample size $n^*$} if the predictive
distributions~$m_i({\mathbf y}^*)$ and~$m_j({\mathbf y}^*)$ are close in terms
of some distance measure for data of that sample size. The model/prior
pairs $\{M_i,\pi_i\}$ and $\{M_j,\pi_j\}$ are \textit{exact predictive
matching at sample size $n^*$} if $m_i({\mathbf y}^*)= m_j({\mathbf y}^*)$ for
all ${\mathbf y}^*$ of sample size~$n^*$.
\end{deft}

One only wants predictive matching for ``minimal'' sample sizes, since,
for larger sample sizes, the discrimination between models occurs
through the marginal densities; they must differ for
discrimination.\vspace*{-1pt}

\begin{crit}[(Predictive matching)]
For appropriately defined ``minimal sample size'' in comparing $M_i$
with $M_j$, one should have model selection priors that are predictive
matching. Optimal (though not always obtainable) is exact predictive matching.\vspace*{-1pt}
\end{crit}

In \citet{BerPer01}, minimal sample size was defined as the smallest
sample size for which the models under consideration have finite
marginal densities when objective estimation priors are used. Typically
this minimal sample size equals the number of parameters in the model
or, more generally, is the number of observations needed for all
parameters to be identifiable.
For model selection, however, minimal sample size needs to be defined
relative to the model selection priors being utilized. Hence we have
the following general definition.\vspace*{-1pt}

%de2 #&#
\begin{deft}[(Minimal training sample)]A \textit{minimal training
sample}~$\mathbf{y}_{i}^*$ for $\{M_i,\pi_i\}$ is a sample of
minimal size
$n_i^*\geq1$ with a finite nonzero marginal density $m_i(\mathbf
{y}_{i}^*)$.\vspace*{-1pt}
\end{deft}

There are many possibilities for even exact predictive matching. We
here highlight two types of exact predictive matching, which are of
particular relevance to the development of objective model selection
priors for the variable selection problem discussed in Section~\ref{priors}.\vspace*{-1pt}

%de3 #&#
\begin{deft}[(Null predictive matching)] \label{defnullpred} The
model selection priors are \textit{null predictive matching} if each of
the model/prior pairs $\{M_i,\pi_i\}$ and $\{M_0,\pi_0\}$ are exact
predictive matching for all minimal training samples $\mathbf
{y}_{i}^*$ for
$\{M_i,\pi_i\}$.\vspace*{-1pt}
\end{deft}
Definition~\ref{defnullpred} reflects the common view---starting with
\citet{Jef61}---that data of a minimal size should not allow one to
distinguish between the null and alternative models. Null predictive
matching arguments have also been used by \citet{GhosSamanta2002} and
\citet{SpSm82} among others.\vspace*{-1pt}
%
%de4 #&#
\begin{deft}[(Dimensional predictive matching)] The model
selection priors are \textit{dimensional predictive matching} if each of
the model/prior pairs $\{M_i,\pi_i\}$ and $\{M_j,\pi_j\}$ of the same
complexity/dimension (i.e., $k_i=k_j$) are exact predictive matching
for all minimal training samples $\mathbf{y}_{i}^*$ for models of
that dimension.\vspace*{-1pt}
\end{deft}
The next section gives the most prominent example of dimensional
predictive matching.\vadjust{\goodbreak}

%s2.5 #&#
\subsection{Invariance criteria}

Invariance arguments have played a prominent role in statistics [cf.
\citet{Ber85}], especially in objective Bayesian estimation theory.
They are also extremely helpful in part of the specification of
objective Bayesian model selection priors.

A basic type of invariance that is almost always relevant for model
selection is invariance to the units of measurement being used:

\begin{crit}[(Measurement invariance)]
The units of measurement used for the observations or model
parameters should not affect Bayesian answers.
\end{crit}

A much more powerful, but special, type of invariance arises when the
family of models under consideration are such that the model structures
are invariant to group transformations. Following the notation in \citet
{Ber85}, we formally state:
%
%de5 #&#
\begin{deft}
The family of densities for $\mathbf{y}\in\mathbb{R}^{n}$,
$\mathfrak{F}:=
\{f(\mathbf{y}\mid\bolds{\theta})\dvtx     \bolds
{\theta}\in\Theta\}$ is said to be {\it
invariant under the group of transformations} ${G}:=\{g\dvtx  \mathbb
{R}^{n}\rightarrow\mathbb{R}^{n}\}$ if, for every $g \in\mathfrak{G}$
and $\bolds{\theta}\in\Theta$, there exists a unique
$\bolds{\theta}^{*} \in
\Theta$ such that $\mathbf{X}=g(\mathbf{Y})$ has density
$f(\mathbf{x}\mid\bolds{\theta}^{*})
\in\mathfrak{F}$. In such a situation, $\bolds{\theta}^{*}$
will be denoted
$\bar{g}(\bolds{\theta})$.
\end{deft}

There are two consequences of applying invariance here. The first is a
new criterion:

\begin{crit}[(Group invariance)]
$\!\!\!$If all models are invariant under a~group of transformations
$G_0$, then the conditional distributions, $\pi_i(\bolds{\beta
}_i\mid\bolds{\alpha}
)$, should be chosen in such a way that the conditional marginal distributions
%
%e5 #&#
\begin{equation}
\label{eqinv-marginals} f_i(\mathbf{y}\mid\bolds{\alpha})=\int
f_i(\mathbf {y}\mid\bolds{\alpha},\bolds{\beta}_i)
\pi_i(\bolds{\beta}_i\mid\bolds{\alpha}) \,d\bolds {
\beta}_i,
\end{equation}
are also invariant under $G_0$. [Here, $(\bolds{\alpha
},\bolds{\beta}_i,i)$ would
correspond to $\bolds{\theta}$ in the definition of invariance.]
\end{crit}

Indeed, the $\pi_i(\bolds{\beta}_i\mid\bolds{\alpha})$
could hardly be called
objective model selection priors if they eliminated an invariance
structure that was possessed by all of the original models. This can
also be viewed as a formalization of
the \citet{Jef61} requirement that the prior for a nonnull parameter
should be ``centered at the simple model.''

The second use of invariance is in determining the objective prior for
the common model parameters $\pi_i(\bolds{\alpha})$. Since all
of the
marginal models, $f_i(\mathbf{y}\mid\bolds{\alpha})$, will
be invariant under $G_0$
if the Group invariance criterion is applied, there are compelling
reasons to choose the prior
%
%e6 #&#
\begin{equation}
\label{eqright-Haar} \pi_i(\bolds{\alpha})= \pi^H(
\bolds{\alpha})\qquad \mbox{for all $i$} ,
\end{equation}
where $\pi^H(\cdot)$ is the right-Haar density corresponding to the
group $G_0$. The reason is given in \citet{BerPerVar98}, namely
that
under commonly satisfied\vadjust{\goodbreak} conditions (satisfied for the variable
selection problem---see Result~\ref{resPM2} in Section~\ref{priors}),
use of a common $\pi^H(\bolds{\alpha})$ for all marginal models
then ensures
exact predictive matching among the models for the minimal training
sample size, as in the example given in Section~\ref{secpmc}.

The most surprising feature of this result is that $\pi^H(\bolds
{\alpha})$
is typically improper (and hence could be multiplied by an arbitrary
constant) and yet, if the same $\pi^H(\bolds{\alpha})$ is used
for all
marginal models, the prior is appropriately calibrated across models in
the strong sense of exact predictive matching. (For any improper prior
that occurred in both the numerator and denominator of a Bayes factor,
any arbitrary multiplicative constant would obviously cancel, but this
is not nearly as compelling a justification as exact predictive
matching.) The right-Haar prior is also the objective estimation prior
for such models, and so has been extensively studied in invariant situations.

Thus, for invariant models, the combination of the Group invariance
criterion and (exact) Predictive matching criterion allows complete
specification of the prior for $\bolds{\alpha}$ in all models.
It is also
surprising that this argument does not require orthogonality of
$\bolds{\alpha}$ and $\bolds{\beta}_i$ (i.e.,
cross-information of zero in the
Fisher information matrix) which, since \citet{Jef61}, has been viewed
as a necessary condition to say that one can use a common prior for
$\bolds{\alpha}$ in different models [see, e.g., \citet
{Hsiao97}, \citet{KaVa92}].

There might be concern here as to use of improper priors, even if they
are exact predictive matching, especially because of the discussion in
Section~\ref{secbasic}. This concern is obviated by the realization
that use of any series of proper priors approximating $\pi^H(\alpha)$
will, in the limit, yield Bayes factors equal
to that obtained directly from $\pi^H(\alpha) $; see Lemma~\ref
{approx-improper} in Appendix~\ref{app1}.

%-------------------------------------Objective prior distributions for
%variable selection in regression models
%s3 #&#
\section{Objective prior distributions for variable selection in normal
linear models}\label{priors}

%s3.1 #&#
\subsection{Introduction}\label{subsecVS}
We now turn to a particular scenario---variable selection in normal
linear models---to illustrate application of the criterion in
Section~\ref{secCriteria}.
Consider a response variable $Y$ known to be explained by $k_0$
variables (e.g., an intercept) and by some subset of $p$ other
possible explanatory variables.
This can formally be stated as a model selection problem with the
following $2^p$ competing models for data $\mathbf{y}=(y_1,\ldots,y_n)$:
%
%e7 #&#
\begin{eqnarray}
\label{eqTheProb}
M_0 \dvtx f_{0}(\mathbf{y}\mid\bolds{
\beta}_{0}, \sigma)&=& \N_n\bigl(\mathbf{y}\mid
\mathbf{X}_0\bolds{\beta}_0,\sigma^2
\mathbf{I}\bigr),
\nonumber
\\[-8pt]
\\[-8pt]
\nonumber
\qquad M_i \dvtx f_{i}(\mathbf{y}\mid\bolds{
\beta}_{i}, \bolds {\beta}_{0}, \sigma)&=&\N_n
\bigl(\mathbf{y}\mid \mathbf{X}_0\bolds{\beta}_0+
\mathbf{X}_i\bolds {\beta}_i,\sigma^2
\mathbf{I}\bigr),\qquad i=1, \ldots, 2^{p}-1 ,
\end{eqnarray}
where $\bolds{\beta}_{0}$, $\sigma$, and the $\bolds
{\beta}_{i}$ are unknown.
Here $\mathbf{X}_0$ is a $n\times k_0$ design matrix corresponding
to the
$k_0$ variables common to all models; often $\mathbf
{X}_0=\mathbf{1}$ so~$M_0$
contains only the intercept. Finally, the $\mathbf{X}_i$ are
$n\times k_i$
design matrices corresponding to $k_i$ of the $p$ other possible
explanatory variables. We make the usual assumption that all design
matrices are full rank (without loss of generality). Note that, if the
covariance matrix is of the form $\sigma^2 \bolds{\Lambda}$
with $\bolds{\Lambda}$ known,\vadjust{\goodbreak}
simply transform $\mathbf{Y}$ so that the covariance matrix is proportional
to the identity;
note that this does not alter the meaning of the $\bolds{\beta
}$'s and hence
the meaning of the models. Also, setting $\bolds{\alpha}=
(\bolds{\beta}_{0},\sigma)$ and $N =2^p$ puts this in the general framework discussed
earlier, with
$M_{0}$ being the null model.

The primary development is for the most common situation of $\sigma$
unknown and $k_0\ge1$, but the simpler cases where either $\sigma$ is
known or $k_0 = 0$ (i.e., the null model only contains the error term)
are briefly treated in Section~\ref{sec2simp}.

%This problem has model space
%$\mathcal{ M} =\{M_{0},\dots, M_{2^{p}-1}\}$ with

In this setting and following Jeffreys desiderata,
\citet{ZellSiow80} recommended use of common objective estimation priors
for $\bolds{\alpha}$ (after orthogonalization) and multivariate Cauchy
priors for $\pi_i(\bolds{\beta}_i\mid\bolds{\alpha})$,
centered at zero and with
prior scale matrix $\sigma^2 n ({\mathbf{X}_i^{\prime
}}{\mathbf{X}_i})^{-1}$; a
similar scale matrix was used in \citet{Zellner86} for the g-prior.

%s3.2 #&#
\subsection{Proposed prior (the ``robust prior'')}

It is useful to first write down the specific form of the prior that
will result
from applying the criteria. Indeed, under model $M_i$, the prior is of
the form
%
%e8 #&#
\begin{eqnarray}
\label{eqourprior}
\pi_i^R(\bolds{\beta}_0,
\bolds{\beta}_i,\sigma)&=&\pi (\bolds{\beta}_0,\sigma)
\times\pi_i^R(\bolds{\beta}_i\mid\bolds{
\beta}_0,\sigma)
\nonumber
\\[-8pt]
\\[-8pt]
\nonumber
&= &\sigma^{-1} \times\int_0^\infty
\N_{k_i}(\bolds{\beta }_i\mid\mathbf{0},g \bolds{
\Sigma}_i) p^{R}_i(g) \,dg,
\end{eqnarray}
where $\bolds{\Sigma}_i=\operatorname{Cov}(\hat{\bolds{\beta }_i})=\sigma^2  (\mathbf{V}_i^t\mathbf{V}_i)^{-1}$ is the
covariance of the maximum likelihood estimator of $\bolds{\beta
}_i$, with
%
%e9 #&#
\begin{equation}
\label{eqVmatrix} \mathbf{V}_{i} = \bigl(\mathbf{I}_{n}-
\mathbf {X}_{0}\bigl(\mathbf{X}_{0}^{t}
\mathbf{X}_{0}\bigr)^{-1}\mathbf {X}_{0}^{t}
\bigr)\mathbf{X}_{i}
\end{equation}
and
%
%e10 #&#
\begin{eqnarray}
\label{eqpiR-g}
&&p_i^{R}(g)=a \bigl[\rho_i(b+n)
\bigr]^{a} (g+b)^{-(a+1)} 1_{\{ g>\rho
_i(b+n)-b \} } ,
\\
\label{eqg-conditions}
&&\qquad\mbox{with } a>0, b>0\quad \mbox{and}\quad \rho_i \geq
\frac{b}{b+n} .
\end{eqnarray}
Note that these conditions ensure that $p_i^{R}(g)$ is a proper
density, and $g$ is positive [necessary in (\ref{eqourprior})], so
that $\pi_i^R(\bolds{\beta}_i\mid\bolds{\beta}_0,\sigma
)$ is proper, satisfying the
first part of the Basic criterion of Section~\ref{secbasic}.
The particular choices of hyperparameters that we favor are discussed
in Section~\ref{sechyperparameters}.

The prior (\ref{eqourprior}) has its origins in the \textit{robust prior}
introduced by \citet{Straw71} and \citeauthor{Ber80} (\citeyear{Ber80,Ber85}), for estimating a
$k$-variate normal mean~$\bolds{\beta}$ in the sampling scheme
$\hat{\bolds{\beta}
}\sim\N_k(\bolds{\beta}, \bolds{\Sigma})$. More
{precisely}, the full conditional
of $\bolds{\beta}_i$ induced by (\ref{eqourprior}) generalizes
the above
mentioned \textit{robust prior} considering the sampling distribution of
the maximum likelihood estimator, namely $\hat{\bolds{\beta
}}_i\sim\N_{k_i}(\bolds{\beta}_i, \sigma^2  (\mathbf
{V}_i^t\mathbf{V}_i)^{-1})$. The primary
reasons for \citet{Straw71} and \citeauthor{Ber80} (\citeyear{Ber80,Ber85}) to consider such
priors was that it results in
closed form inferences, including closed form Bayes factors, and
results in estimates that are robust in various senses. For this
reason, we continue the tradition of calling (\ref{eqourprior}) the
\textit{robust prior} and use a superindex~$R$ to denote it. Note also
that priors of this form have been previously\vadjust{\goodbreak} considered.
The priors proposed by \citet{liang08} are particular cases with
$a=1/2,  b=1,  \rho_i=1/(1+n)$ (the hyper-g prior) and $a=1/2,
b=n,  \rho_i=1/2$ (the hyper-g/n prior). The prior in \citet{CuiGeo08}
has $a=1,  b=1, \allowbreak \rho_i=1/(1+n)$. The original Berger prior for
robust estimation is the particular case with $a=1/2,  b=1,  \rho_i=(k_i+1)/(k_i+3)$; closely related priors are those of \citet
{MarStr10},
\citet{MarGeo08}.

Finally, it is useful to note that $\pi_i^R(\bolds{\beta}_i\mid
\bolds{\beta}_0,\sigma)$ behaves in the tails as a~multivariate Student distribution
(already noticed for a particular case in \citet{Ber80}, and the reason
for its robust estimation properties).

%pr1 #&#
\begin{prop}\label{piRSt} Writing $\|\bolds{\beta}_i\|^2=\bolds{\beta}_i^t(\mathbf{V}_i^t\mathbf
{V}_i)\bolds{\beta}_i$,
\[
\lim_{\|\beta_i\|^2\rightarrow\infty} \frac{\pi_i^R(\bolds{\beta}\mid\bolds{\beta}_0,\sigma)}{{\operatorname
St}_{k_i}(\bolds{\beta}\mid\mathbf{0},   (a \Gamma
(a))^{1/a}   \rho_i   \mathbf{B}^*(b, \sigma)/a,   2a)} = 1,
\]
where $\mathbf{B}^*(b, \sigma)=\sigma^2 (b+n) (\mathbf
{V}_i^t\mathbf{V}_i)^{-1}$.
\end{prop}
\begin{pf} See Appendix~\ref{app2}.
\end{pf}

In the model selection scenario, the thickness of the prior tails is
related to the information consistency criteria, and is the reason
\citet{Jef61} used a Cauchy as the prior for testing a normal mean.
Also, using this result, we can see that $\pi_i^R(\bolds{\beta
}_i\mid\bolds{\beta}_0,\sigma)$ has close connections with the Zellner--Siow priors; in
fact, for $a=1/2$, $b=n$, $\rho_i=2/\pi$ and large $n$, $\pi_i^R(\bolds{\beta}_i\mid\bolds{\beta}_0,\sigma)$, and the Zellner--Siow priors
have exactly the
same tails.

%--------------------------------------------------------------
%s3.3 #&#
\subsection{\texorpdfstring{Justification of model selection priors of the form (\protect\ref{eqourprior})}
{Justification of model selection priors of the form (8)}}

We will use the Group invariance criterion and Predictive matching
criterion (along with practical computational considerations)
to justify use of model selection priors of the form~(\ref
{eqourprior}). We first justify the use of
$\pi^R(\bolds{\beta}_{0},\sigma) =1/\sigma$ for the common
parameters and then
justify the choice $\pi_i^R(\bolds{\beta}\mid\bolds
{\beta}_0,\sigma)$ for the
model specific parameters.

%s3.3.1 #&#
\subsubsection{Justification of the prior for the common parameters}
It is convenient, in this section, to consider a more general class of
conditional priors,
%
%e12 #&#
\begin{equation}
\label{eqscale} \pi_i(\bolds{\beta}_{i}\mid\bolds{
\beta}_{0}, \sigma )=\sigma^{-k_{i}}h_{i}\biggl(
\frac
{\bolds{\beta}_{i}}{\sigma}\biggr) ,
\end{equation}
where $h_i$ is any proper density with support $\mathcal{ R}^{k_i}$.
The robust prior is the particular case
%
%e13 #&#
\begin{equation}
\label{eqmixture} h_{i}^R(\mathbf{u})= \int
\N_{k_{i}}\bigl(\mathbf{u}\mid\mathbf{0}, g \bigl(\mathbf
{V}_i^t\mathbf{V}_i\bigr)^{-1}
\bigr) p_{i}^R(g) \,dg .
\end{equation}

It is shown, in Appendix~\ref{app3}, that all models in (\ref{eqTheProb}) are
invariant under the group of transformations
\[
G_0=\bigl\{g=(c,\mathbf{b}) \in(0,\infty)\times\mathcal{
R}^{k_0}\dvtx g(\mathbf{y})\rightarrow c\mathbf{y}+\mathbf{X}_0
\mathbf{b}\bigr\}.
\]
The following establishes a necessary and sufficient condition on the
conditional prior $\pi_i(\bolds{\beta}_i\mid\bolds{\beta
}_{0},\sigma)$ for the
Group invariance criterion to hold for this group.\vspace*{-1pt}

%re1 #&#
\begin{res}\label{resIcrit}
The conditional marginals
%
%e14 #&#
\begin{equation}
\label{eqMiI} f_i(\mathbf{y}\mid\bolds{\beta}_{0},
\sigma)=\int\N_n\bigl(\mathbf{y}\mid\mathbf{X}_0\bolds{
\beta}_{0} +\mathbf{X}_i\bolds{\beta}_i,
\sigma^2\mathbf {I}\bigr) \pi_i(\bolds{
\beta}_i\mid\bolds{\beta}_{0}, \sigma) \,d\bolds{
\beta}_i
\end{equation}
are invariant under $G_0$ if and only if $\pi_i(\bolds{\beta
}_i\mid\bolds{\beta}_{0}, \sigma)$ has the form \eqref{eqscale}.\vspace*{-1pt}
\end{res}
\begin{pf} See Appendix~\ref{app3}.\vspace*{-1pt}
\end{pf}

Based on the Group invariance criterion, Result~\ref{resIcrit} implies
that, conditionally on the common parameters $\bolds{\beta}_0$
and $\sigma$,
$\bolds{\beta}_i$ must be scaled by $\sigma$, centered at zero
and not depend
on $\bolds{\beta}_0$ [as was argued for simple normal testing in
\citet
{Jef61}]. Note, in particular, that the robust prior in (\ref
{eqourprior}) satisfies the Group invariance criterion (although it is
not the only prior that does so).

Next, since each marginal model $f_i(\mathbf{y}\mid\bolds
{\beta}_{0},  \sigma)$
resulting from a prior in~\eqref{eqscale} is invariant with respect to
$G_0$, the suggestion from \citet{BerPerVar98} is to use the right-Haar
density for the common parameters $(\bolds{\beta}_{0},\sigma)$, namely
\[
\pi_i(\bolds{\beta}_{0},\sigma)=\pi^H(
\bolds{\beta }_{0},\sigma)=\sigma^{-1},
\]
the right-Haar prior for the location-scale group. Using this, the
overall model prior would be of the form
%
%e15 #&#
\begin{equation}
\label{eqscale2} \pi_i(\bolds{\beta}_0,\bolds{
\beta}_i,\sigma)=\sigma^{-1-k_i} h_i\biggl(
\frac{\bolds{\beta}_i}{\sigma}\biggr).
\end{equation}
The justification for the right-Haar prior in \citet{BerPerVar98}
depends, however,
on showing that it is predictive matching, in the sense described in
the following result.\vspace*{-1pt}

%re2 #&#
\begin{res}\label{resPM2}
For $M_i$, let the prior $\pi_i(\bolds{\beta}_0,\bolds
{\beta}_i,\sigma)$ be of the
form \eqref{eqscale2}, where~$h_i$ is symmetric about zero. Then all
model/prior pairs $\{M_i,\pi_i\}$ are exact predictive matching for
$n^{*}=k_{0}+1$.\vspace*{-1pt}
\end{res}

\begin{pf} See Appendix~\ref{app4}.\vspace*{-1pt}
\end{pf}

The conclusion of the above development is that the Group invariance
criterion and Predictive matching criterion imply that model selection
priors should be of the form \eqref{eqscale2}, with $h_i$ symmetric
about zero. It would thus appear that the robust prior satisfies these
criteria, as (\ref{eqmixture}) is clearly symmetric about zero. [Any\vadjust{\goodbreak}
scale mixture of Normals would also satisfy these criteria, since the
resulting $h(\cdot)$ would be symmetric about 0.]
Note, however, that~$h_i^R$ has scale matrix proportional to
$(\mathbf{V}_i^t\mathbf{V}_i)^{-1}$, and $\mathbf{V}_i$
in (\ref{eqVmatrix}) requires both~$\mathbf{X}_0$ and $\mathbf{X}_i$, which would seem to
indicate that a sample size
of $k_0+k_i$ is required. Hence, Result~\ref{resPM2} would seem to
apply to the robust prior only if $k_i=1$.

This is a situation, however, where the definition of a minimal sample
size is somewhat ambiguous. For instance, suppose one were presented
$\mathbf{X}_0$ and $\mathbf{X}_i$ for $k_0+k_i$ observations
for each model $M_i$,
but that only $k_0+1$ of the $y_i$ was reported for all models, with
the rest being missing data. This is still a minimal sample size in the
sense that it is the smallest collection of $y_i$ for which all
marginal densities exist for the robust prior, and now Result \ref
{resPM2} applies to say that the robust prior is predictive matching
for all models.

%s3.3.2 #&#
\subsubsection{Justification of the prior for the model specific parameters}
While the robust prior is thus validated as satisfying the group
invariance criterion and a version of the predictive matching
criterion, there are many other model selection priors of form (\ref
{eqscale2}) which also satisfy these criteria. There are additional
reasons, however, to focus on the robust priors with $h^R_i(\mathbf
{u})$ of
form (\ref{eqmixture}). The first is that only scale mixtures of
normals seem to have any possibility of yielding Bayes factors that
have closed form. While we have not focused on this as a necessary
criterion, it is an attractive enough property to justify the
restriction. There are, however, two other features of~(\ref
{eqmixture}) that need justification: the use of the mixture density $
p^R_{i}(g)$, and the choice of the conditional scale matrix
$(\mathbf{V}_i^t\mathbf{V}_i)^{-1}$.

The mixture density $p^R_{i}(g)$ encompasses virtually all of the
mixtures that have been found which can lead to closed form expressions
for Bayes factors; for example, Zellner--Siow priors are scale mixtures
of normals, but with a~different mixing density which does not lead to
close-form expressions. [The choice of mixing density in \citet
{MarGeo08} is a~very interesting exception, in that it leads to a~closed
form expression for a~different reason than does $p^R_{i}(g)$.]
So, while not completely definitive, $p^R_{i}(g)$ is an attractive
choice. The choice of $(\mathbf{V}_i^t\mathbf{V}_i)^{-1}$ as
the conditional scale
matrix seems much more arbitrary, but there is one standard argument
and one surprising argument in its favor.

The standard argument is the measurement invariance criterion; if the
conditional scale matrix is chosen to be $(\mathbf
{V}_i^t\mathbf{V}_i)^{-1}$, it is
easy to see that Bayes factors will be unaffected by changes in the
units of measurement of either $\mathbf{y}$ or the model
parameters. But
there are many other choices of the conditional scale matrix which also
have this property.

A quite surprising predictive matching result that supports use of
$(\mathbf{V}_i^t\mathbf{V}_i)^{-1}$ as the conditional scale
matrix is as follows.
%
%re3 #&#
\begin{res}\label{resPM1}
For $M_i$, let the prior be as in \eqref{eqscale2}, where $h_i$ is the
scale mixture of normals in \eqref{eqmixture}. The priors are
then\vadjust{\goodbreak}
null predictive matching and dimensional predictive matching for
samples of size $k_0+k_i$, and no choice of the conditional scale
matrix other than $(\mathbf{V}_i^t\mathbf{V}_i)^{-1}$ (or a
multiple) can achieve
this predictive matching.
\end{res}

\begin{pf}See Appendix~\ref{app5}.
\end{pf}

This is surprising, in that it is a predictive matching result for
larger sample sizes ($k_0+k_i$) than are encountered in typical
predictive matching results, such as Result~\ref{resPM2}. That it only
holds for conditional scale matrices proportional to $(\mathbf
{V}_i^t\mathbf{V}_i)^{-1}$ is also surprising, but does strongly
support choosing a
prior of the form \eqref{eqourprior}.

%s3.4 #&#
\subsection{Choosing the hyperparameters for $p^R_{i}(g)$}
\label{sechyperparameters}

%s3.4.1 #&#
\subsubsection{Introduction}
The Bayes factor of $M_i$ to $M_0$ arising from the robust prior $\pi_i^R$ in \eqref{eqourprior} can be compactly expressed as the
following function of the hyperparameters $a$, $b$ and $\rho_i$:
%
%e16 #&#
\begin{equation}
\label{eqBFR} B_{i0}=Q_{i0}^{-{(n-k_0)}/{2}} \frac{2a}{k_i+2a}
\bigl[ \rho_i (n + b) \bigr]^{-k_i/2} \operatorname{AP}_i,
\end{equation}
where $\operatorname{AP}_i$ is the hypergeometric function of two variables [see
\citet{Weis}], or Apell hypergeometric function
\[
\label{eqapell}
\operatorname{AP}_i =\mathrm{F}_1
\biggl[a+\frac{k_i}{2}; \frac
{k_i+k_0-n}{2}, \frac{n-k_0}{2}; a+1+
\frac{k_i}{2}; \frac{(b-1)}{\rho_i(b+n)}; \frac
{b-Q_{i0}^{-1}}{\rho_i(b+n)} \biggr] ,
\]
and $Q_{i0}=\operatorname{SSE}_i/\operatorname{SSE}_0$ is the ratio of the sum of squared errors of
models~$M_i$ and~$M_0$. The details of this computation are given in
Appendix~\ref{app6}.

Having a closed form expression for Bayes factors is not one of our
formal criteria for model selection priors, but it is certainly a
desirable property, especially when realizing that one is dealing with
$2^p$ models in variable selection.

The values for the hyperparameters that will be recommended are
$a=1/2$, $b=1$ and $\rho_i=(k_i+k_0)^{-1}$. The arguments justifying
this specific recommendation follow.

%s3.4.2 #&#
\subsubsection{Implications of the consistency criteria}
The consistency criteria of Section~\ref{secfcrit} provide
considerable guidance as to the choice of $a$, $b$ and the~$\rho_i$.
In particular, they lead to the following result.

%re4 #&#
\begin{res}\label{resconsistencycond}
The three consistency criterion of Section~\ref{secconsiscrit} are
satisfied by the robust prior if
$a$ and $\rho_i$ do not depend on $n$, $\lim_{n\rightarrow\infty
}\frac
{b}{n} = c \geq0$, $\lim_{n\rightarrow\infty}\rho_{i} (b+n) =
\infty
$ and $n\geq k_i+k_0+2a$.
\end{res}

This result follows from (\ref{eqmodconsist}), (\ref
{eqIPconditions}) and (\ref{eqinfconsist}) below, which are presented
as separate results because they can be established in more generality
than simply for the robust prior.

%pa3.4.2.1 #&#
\paragraph*{Use of model selection consistency}

Suppose $M_i$ is the true model, and consider any other model $M_j$.
A key assumption for model selection consistency [\citet{FLS01}] is
that, asymptotically, the design matrices are such that the models are
differentiated, in the sense that
%
%e17 #&#
\begin{equation}
\label{eqdifferentiate} \lim_{n \rightarrow\infty} \frac{\bolds{\beta}_i^{t}
\mathbf{V}_i^{t} (\mathbf{I} - \mathbf{P}_j) \mathbf
{V}_i \bolds{\beta}_i}{n} =
b_j \in(0,\infty) ,
\end{equation}
where $\mathbf{P}_j = \mathbf{V}_j (\mathbf
{V}_j^{t}\mathbf{V}_j)^{-1} \mathbf{V}_j^{t}$.
%
%re5 #&#
\begin{res}\label{resCM1}
Suppose (\ref{eqdifferentiate}) is satisfied and that the priors $\pi_i(\bolds{\beta}_0,\bolds{\beta}_i,\sigma)$ are of the
form \eqref{eqscale2}, with
$h_{i}(\mathbf{u})=
\int\N_{k_{i}}(\mathbf{u}\mid\mathbf{0}, g (\mathbf
{V}_i^t\mathbf{V}_i)^{-1})  p_{i}(g)
\,dg$. If the
$p_{i}(g)$ are proper densities such that
\[
\lim_{n \rightarrow\infty} \int_0^\infty(1+g)^{-{k_i}/{2}}
p_i(g) \,dg =0 ,
\]
model selection consistency will result.
\end{res}

\begin{pf} The proof follows directly from the proof of Theorem 3 in
\citet{liang08} and is, hence, omitted.
\end{pf}

%co1 #&#
\begin{cor}\label{resCM2}
The prior distributions in \eqref{eqourprior} are model selection
consistent if
%
%e18 #&#
\begin{equation}
\label{eqmodconsist} \lim_{n\rightarrow\infty}\rho_{i} (b+n)=\infty .
\end{equation}

\end{cor}

\begin{pf} See Appendix~\ref{app7}.
\end{pf}

%pa3.4.2.2 #&#
\paragraph*{Use of intrinsic prior consistency}
Related to (\ref{eqdifferentiate}) is the condition that
%
%e19 #&#
\begin{equation}
\label{intrinsic-prior} \lim_{n \rightarrow\infty} \frac{1}{n}
\mathbf{V}_l^{t}\mathbf{V}_l = \bolds{
\Xi}_l
\end{equation}
for some positive definite matrix $\bolds{\Xi}_l$. This would trivially
happen if either there is a fixed design with replicates, or when the
covariates arise randomly from a fixed distribution having second moments.

%re6 #&#
\begin{res}\label{resIP}
If (\ref{intrinsic-prior}) holds,
%
%e20 #&#
\begin{equation}
\label{eqIPconditions} \mbox{$a$ and $\rho_i$ do not depend on
$n$}\quad \mbox{and}\quad \frac{b}{n} \rightarrow c ,
\end{equation}
then the conditional robust prior $\pi_i^R(\bolds{\beta}_i\mid
\bolds{\beta}_0,\sigma
)$ in (\ref{eqourprior}) converges to the fixed intrinsic prior
%
%e21 #&#
\begin{equation}
\label{eqIP-normal} \pi_i(\bolds{\beta}_i\mid\bolds{
\beta}_0,\sigma)=\int_0^\infty
\N_{k_i}\bigl(\bolds{\beta}_i\mid\mathbf{0},g^*
\sigma^2 \bolds{\Xi}^{-1}\bigr) p_i\bigl(g^*
\bigr) \,dg^*,
\end{equation}
where $p_i(g^*)=a[\rho_i(c+1)]^{a}  (g^*+c)^{-(a+1)} 1_{\{ g^*>\rho
_i(c+1)-c\}}$.\vadjust{\goodbreak}
\end{res}

\begin{pf} Changing variables to $g^*=g/n$, the integral in (\ref
{eqourprior}) becomes
\begin{eqnarray*}
&&\int_0^\infty\N_{k_i} \biggl(\bolds{
\beta}_i \Bigm|\mathbf {0},g^*\sigma^2 \biggl(
\frac{1}{n} \mathbf{V}_l^{t}\mathbf{V}_l
\biggr)^{-1} \biggr) a \biggl[\rho_i \biggl(
\frac{b}{n}+1 \biggr) \biggr]^{a}
\\
&&\qquad{}\times \biggl(g^*+\frac{b}{n} \biggr)^{-(a+1)} 1_{\{ g^*>\rho_i(
{b}/{n}+1)-{b}/{n}\}}
\,dg^* .
\end{eqnarray*}
For large $n$ and using (\ref{intrinsic-prior}) and (\ref
{eqIPconditions}), it is easy to find an integrable function
dominating the integrand, so the dominated
convergence theorem can be applied to interchange the integral and
limit, yielding the result.
\end{pf}

%pa3.4.2.3 #&#
\paragraph*{Use of information consistency}
For the variable selection problem, it is easy to see that
\[
\sup_{\bolds{\beta}_{l}, \bolds{\beta}_{0}, \sigma
}f_{l}(\mathbf{y}\mid\bolds{\beta}_{0},
\bolds {\beta}_{l}, \sigma)=(2\pi\operatorname{SSE}_{l}/n)^{-n/2}
\exp (-n/2 )
\]
for model $M_l$. Hence, for any given data set $\mathbf{y}$, the estimated
likelihood ratio in (\ref{eqlr}) is
\[
\Lambda_{i0}(\mathbf{y})=Q_{i0}(\mathbf{y})^{-n/2} ,
\]
where $Q_{i0}(\mathbf{y})$ is the ratio of the residual sum of
squares of the
two models
for $\mathbf{y}$. %$If $M_{j}$ is nested within $M_{i}$ (so that
%$k_{i}>k_{j}$),
{ Therefore, having a sequence of data vectors $\{\mathbf{y}_{m}\}
$ such that
$\lim_{m\rightarrow\infty}\Lambda_{i0}(\mathbf{y}_{m})= \infty
$ is equivalent
to having a sequence of data vectors such that $\lim_{m\rightarrow
\infty}Q_{i0}(\mathbf{y}_{m})\rightarrow0$.}
%
%re7 #&#
\begin{res}\label{resCI} If $\rho_i \geq b/(b+n)$, the prior in
\eqref{eqourprior}
results in an information consistent Bayes factor for $M_{i}$ versus
%%the nested model $M_{j}$,
{ $M_0$,} if and only if
%
%e22 #&#
\begin{equation}
\label{eqinfconsist} n\geq k_i+k_0+2a .
\end{equation}
\end{res}
\begin{pf} See Appendix~\ref{app8}.
\end{pf}

%s3.4.3 #&#
\subsubsection{Specific choices of hyperparameters}

%pa3.4.3.1 #&#
\paragraph*{The choice of $a$}
Note that, with $k_i > k_j$ and $n \geq k_i+k_0+1$, the Bayes factor
$B_{ij}$ between $M_i$ and $M_j$ exists. It is desirable to
have information consistency for all such sample sizes, in which case
(\ref{eqinfconsist}) would require $ a \leq1/2$.
The choice $a=1/2$ is attractive, in that it coincides with the choice
in \citet{Ber85} and, with this choice, $\pi_l^R$ has Cauchy tails, as
do the popular proposals of \citet{Jef61} and \citeauthor{ZellSiow80} (\citeyear{ZellSiow80,ZellSiow84}).

Additional motivation for this choice can be found by studying the
behavior of $B_{i0}$ when the information favors $M_0$, in the sense
that $Q_{i0} \rightarrow1$. Indeed, \citet{Forte11} shows that the
limiting value of $B_{i0}$ is then bounded above by $2a/(2a+k_i)$ for
any sample size, including a small sample size such as $k_0+k_i+1$. A
small value of $a$ would imply strong evidence in favor of $M_0$, which
does not seem reasonable\vadjust{\goodbreak} when the sample size is small. In contrast,
the recommended choice would yield a bound of $1/(1+k_i)$, which
certainly favors $M_0$, but in a sensibly modest fashion when the
sample size is small.

%pa3.4.3.2 #&#
\paragraph*{The choice of $b$}
To understand the effect of $b$ and the $\rho_i$ on the robust prior,
it is useful to begin by considering the approximating intrinsic prior
in Result~\ref{resIP}, which depends on the hyperparameters only
through the mixing distribution $ p^R_i(g^*)$, which for $a=1/2$ is
given by (when $b/n \rightarrow c$)
%
%e23 #&#
\begin{equation}
\label{eqIPg} { p^R_i\bigl(g^*\bigr)}=
\tfrac{1}{2}\bigl[\rho_i(c+1)\bigr]^{1/2} \bigl(g^*+c
\bigr)^{-3/2} 1_{\{ g^*>\rho_i(c+1)-c\}} .
\end{equation}
This is a very flat-tailed distribution with median $4\rho_i(1+c)-c$.
Because it is so flat tailed, the choice of $c$ in $(g^*+c)^{-3/2}$ is
not particularly influential, so that the main issue is the choice of
the median. For selecting a~median, however, $\rho_i$ and~$c$ are
confounded; that is, we do not need both. For simplicity, therefore, we
will choose $c=0$ (i.e., $b$ such that $b/n \rightarrow0$).

If $b/n \rightarrow c=0$, the intrinsic prior does not depend at all on
$b$. Furthermore, there is very little dependence on $b$, in this case,
for the actual robust prior, as was verified for moderate and small $n$
in \citet{Forte11} through an extensive numerical study.

Since any choice of $b$ for which $b/n \rightarrow0$ makes little
difference, it would be reasonable to make such a choice based on
pragmatic considerations. In this regard, note that the choice $b=1$
has a notable computational advantage, in that the hypergeometric
function of two variables, $\operatorname{AP}_i$, then becomes the standard
hypergeometric function of one variable [\citet{AS64}]. We thus choose $b=1$.

%pa3.4.3.3 #&#
\paragraph*{The choice of $\rho_i$}
This is the most difficult choice to make, since there is only limited
guidance from the various criteria. To review (and assuming $b=1$), we
have that $\rho_i \geq1/(1+n)$ (so that $g>0$); $\lim_{n\rightarrow
\infty} \rho_{i} (1+n)=\infty$ (for model selection consistency); and
$\rho_i$ should not depend on $n$ (for there to be a limiting intrinsic
prior). Also note that $n$ is necessarily greater than or equal to
$k_0+k_i$ for the robust prior and marginal likelihood to exist;
supposing we wish to choose $\rho_i$ so that the conditions are
satisfied for all such $n$, these restrictions only imply that
\[
\mbox{$\rho_i$ must be a constant (independent of $n$) and $
\rho_i \geq 1/(1+k_0+k_i)$} .
\]
We present two arguments below for the specific choice $\rho_i = 1/(k_0+k_i)$.

\begin{argg} Consider the Bayes factor $B_{i0}$ of $M_i$ to $M_0$.
In Result~\ref{resPM1}, it was established that $B_{i0}=1$ for a
sample of size $n=k_{i}+k_{0}$, but a natural question is---what
should we expect for a sample of size $n=k_{i}+k_{0}+1$? Can a single
additional observation provide much information to discriminate between
$M_{i}$ and $M_{0}$? Intuition says no. To quantify the intuition,
consider the situation in which $Q_{i0}\rightarrow1$, which
corresponds to information being as supportive as possible of $M_0$. It
is straightforward to show that,\vadjust{\goodbreak} when $n=k_{i}+k_{0}+1$,
%
%e24 #&#
\begin{equation}
\label{eqrhobound} \lim_{Q_{i0}\rightarrow1} B_{i0}^R =
\frac{1}{k_i+1}\bigl[\rho_i(k_{i}+k_{0}+2)
\bigr]^{-k_i/2} .
\end{equation}
As we should not expect a single extra observation to provide very
strong evidence, even in the case that $Q_{i0}\rightarrow1$, the
implication is that we should choose $\rho_i$ to be as small as is
reasonable. The choice $\rho_i=1/(k_0+k_i+1)$ is the minimum value of~$\rho_i$ and is, hence, certainly a candidate.
\end{argg}

\begin{argg} Consider the intrinsic prior defined by (\ref
{eqIP-normal}) and (\ref{eqIPg}). Note that we have chosen $c=0$
(through the choice of $b=1$) and, after making the transformation
$\tilde g = g^*/\rho_i$, the intrinsic prior can be written
%
%e25 #&#
\begin{equation}
\label{eqtransf-intrinsic} \pi_i(\bolds{\beta}_0,
\bolds{\beta}_i,\sigma)=\sigma^{-1}\times\int
_0^\infty\N_{k_i}\bigl(\bolds{
\beta}_i\mid\mathbf{0},{\tilde g} \rho_i \bolds{
\Xi}^{-1}\bigr) p_i(\tilde g) \,d\tilde g ,
\end{equation}
where $p_i(\tilde g)=(1/2) (\tilde g)^{-3/2} 1_{\{\tilde g>1\}}$.
Thus we see that, in the intrinsic prior approximation to the robust
prior, $\rho_i$ can be interpreted as simply a scale factor to the
conditional covariance matrix. This helps, in that there have been
previous suggestions related to ``unit information priors'' [\citet
{KassWass95}, \citet{BerBarPer10}]. For instance, \citet{BerBarPer10} consider
the group means problem defined as follows: the observations are
\[
y_{ij} = \mu_i + \varepsilon_{ij} ,\qquad i=1,
\ldots, k \mbox{ and }  j=1, \ldots, r ,
\]
with i.i.d. $\varepsilon_{ij}\sim N(\cdot\mid0, \sigma^2)$. Thus there
are $k$ different means, $\mu_i$, and $r$ replicate observations for
each. Applying the robust prior to this example (considering the full
model with all $\mu_i$) results in a conditional covariance matrix in
(\ref{eqtransf-intrinsic}) of $\rho  k   \mathbf{I}$, which
is much too
diffuse if $k$ is large and $\rho$ is not small. Selecting $\rho=1/k$,
on the other hand, restores a ``unit information'' prior. Here $k_0=0$,
so the choice $\rho=1/k$ is equivalent to the overall choice $\rho_i
=1/(k_0+k_i)$. This overall choice is obviously very close to earlier
suggested $1/(k_0+k_i+1)$.
\end{argg}

%------------------------------- Two Simpler
%Cases-------------------------------------
%s3.5 #&#
\subsection{Two simpler cases}\label{sec2simp}

We conclude with discussion of the modifications of the robust prior
that are needed when $\bolds{\beta}_{0}=0$ or when $\sigma$ is known.
%---------------------------------------
%s3.5.1 #&#
\subsubsection{\texorpdfstring{When $\bolds{\beta}_{0}=\mathbf{0}$ and
$\sigma$ is unknown}{When beta0 = 0 and
sigma is unknown}}
When $\bolds{\beta}_{0}=\mathbf{0}$, the robust prior
distribution is
\[
\pi_i^R(\bolds{\beta}_i,\sigma)=\pi(
\sigma) \times\pi_i^R(\bolds{\beta}_i\mid
\sigma)=\sigma^{-1} \times\int_0^\infty
\N_{k_i}(\bolds {\beta}_i\mid\mathbf{0},g \bolds{
\Sigma}_i) p^{R}_i(g) \,dg,
\]
where $\bolds{\Sigma}_i=\operatorname{Cov}(\hat{\bolds{\beta
}_i})=\sigma^2  (\mathbf{X}_i^t\mathbf{X}_i)^{-1}$, the
covariance of the maximum likelihood estimator of $\bolds{\beta
}_i$ and, as before,
\[
p_i^{R}(g)=a\bigl[\rho_i(b+n)
\bigr]^{a} (g+b)^{-(a+1)}, \qquad g>\rho_i(b+n)-b.
\]
The corresponding Bayes factor is as in \eqref{eqBFR} with $k_{0}=0$;
when we choose $a=1/2$, $b=1$ and $\rho_{i}=1/(k_{i}+k_{0})$, it
assumes the simpler form in \eqref{eqRBayesFactorR}, again with $k_{0}=0$.

In regards to the group invariance criterion, when $\bolds{\beta
}_{0}=0$ the
models are invariant under the scale group of transformations,
$G_0=\{\mathbf{y}\rightarrow c\mathbf{y}, c>0\}
$, and it is easy to
show that $\pi(\bolds{\beta}_{i}\mid\bolds{\beta}_{0},
\sigma)$ still needs to be
a scale prior, as in \eqref{eqscale}, to preserve the invariance
structure; also, the use of $\pi(\sigma)=1/\sigma$ is again justified
by predictive matching, as it is the Haar prior for the group.
Null and dimensional predictive matching also hold as well as the
various consistency criteria.

%------------------------------------------
%s3.5.2 #&#
\subsubsection{\texorpdfstring{When $\sigma$ is known and $\bolds{\beta}_{0}
\neq0$}{When sigma is known and beta0 /= 0}}
When $\sigma$ is known, the robust prior becomes
\[
\pi_i^R(\bolds{\beta}_i, \bolds{
\beta}_{0}, \sigma )=\pi(\bolds{\beta}_{0}) \times
\pi_i^R(\bolds{\beta}_i\mid\bolds{
\beta}_{0})\propto\int_0^\infty
\N_{k_i}(\bolds{\beta}_i\mid\mathbf{0},g \bolds{
\Sigma}_i) p^{R}_i(g) \,dg,
\]
where $\bolds{\Sigma}_i=\operatorname{Cov}(\hat{\bolds{\beta
}_i})=\sigma^2  (\mathbf{V}_i^t\mathbf{V}_i)^{-1}$, and
$p_i^{R}(g)$ is as before.

The models are now invariant under the location group $G_0=\{
\mathbf{y}\rightarrow\mathbf{y}+\mathbf{X}_{0}
\mathbf{b},  \mathbf{b}\in\mathcal{
R}^{k_0}\}$,
and it is easy to show that $\pi(\bolds{\beta}_{i}\mid
\bolds{\beta}_{0})$ just
needs to be independent of $\bolds{\beta}_{0}$ to preserve the invariance
structure; the use of the Haar prior $\pi(\bolds{\beta}_{0})=1$
is again
justified through predictive matching arguments.

The Bayes factor can be expressed as
\[
B_{i0}=\int_{0}^{\infty}(g+1)^{-k_{i}/2}
\Lambda_{0i}^{ (
{1}/{(g+1)}-1 )}p_{i}(g)\,dg ,
\]
where $\Lambda_{0i}= \exp (-[\operatorname{SSE}_{0}-\operatorname{SSE}_{i}]/(2\sigma^{2})
)$. This is curiously difficult to express in closed-form in general
but, for our preferred choice $b=1$, change of variables to $h=1/(1+g)$ yields
\begin{eqnarray*}
B_{i0} &=&\int_{0}^{\infty}(g+1)^{-k_{i}/2}
\Lambda_{0i}^{ ({1}/{(g+1)}-1 )}a\bigl(\rho_i(1+n)
\bigr)^{a} (g+1)^{-(a+1)} 1_{\{g>\rho_i(1+n)-1\}}\,dg
\\
&=&a\bigl(\rho_i(1+n)\bigr)^{a}\Lambda_{0i}^{-1}
\int_{0}^{1/[\rho_i(1+n)]}h^{(a-1
+k_{i}/2)}e^{-h[\operatorname{SSE}_{0}-\operatorname{SSE}_{i}]/(2\sigma^{2})}\,dh
\\
&=& a\bigl(\rho_i(1+n)\bigr)^{a}\Lambda_{0i}^{-1}
\biggl(\frac
{[\operatorname{SSE}_{0}-\operatorname{SSE}_{i}]}{2\sigma^{2}} \biggr)^{-(a-2+{k_i}/{2})}
\\
&&{} \times \biggl(\Gamma \biggl[a+\frac{k_i}{2} \biggr]- \Gamma \biggl[a+
\frac{k_i}{2},\frac{[\operatorname{SSE}_{0}-\operatorname{SSE}_{i}]}{2\sigma^{2}\rho_i(1+n)} \biggr] \biggr),
\end{eqnarray*}
where $\Gamma(\nu_1, \nu_2)$ is the incomplete gamma function,
\[
\Gamma(\nu_1, \nu_2)=\int_{\nu_2}^\infty
t^{\nu_1-1}e^{-t} \,dt .
\]
All of the properties of the procedures for the $\sigma$ unknown case
also hold here, except for null predictive matching.

%s4 #&#
\section{Methodological summary for variable selection}
\label{secsummary}

Although the primary purpose of the paper was to develop the criteria
for choice of model selection priors and study their implementation in
an example, the methodological results obtained for the problem of
variable selection in the normal linear model, as outlined in
Section~\ref{subsecVS}, are of interest in their own right.
For ease of use, we summarize these developments here.

Using the notation of Section~\ref{subsecVS}, the prior distribution
recommended for the parameters under model $M_i$ is
\[
\pi_i^R(\bolds{\beta}_0,\bolds{
\beta}_i,\sigma)=\sigma^{-1} \times\int
_0^\infty\N_{k_i}(\bolds{
\beta}_i\mid\mathbf{0},g \bolds {\Sigma}_i)
p^{R}_i(g) \,dg,
\]
where $\bolds{\Sigma}_i=\sigma^2  (\mathbf
{V}_i^t\mathbf{V}_i)^{-1}$, $\mathbf{V}_{i} = (\mathbf
{I}_{n}-\mathbf{X}_{0}(\mathbf{X}_{0}^{t}\mathbf
{X}_{0})^{-1}\mathbf{X}_{0}^{t})\mathbf{X}_{i}$,
and
%
%e26 #&#
\[
p_i^{R}(g)=\frac{1}{2} \biggl[
\frac{(1+n)}{(k_i+k_0)} \biggr]^{
{1}/{2}} (g+1)^{-3/2} 1_{\{ g>(k_i+k_0)^{-1}(1+n)-1 \} } .
\]

The resulting Bayes factors have closed form expressions in terms of
the the hypergeometric function, namely
%
%e27 #&#
\begin{eqnarray}
\label{eqRBayesFactorR} \qquad B_{i0}&=& \biggl[ \frac{n + 1}{k_i +k_0}
\biggr]^{-{k_i}/{2}}
\nonumber
\\[-8pt]
\\[-8pt]
\nonumber
&&{}\times\frac
{Q_{i0}^{-{(n-k_0)}/{2}}}{k_i+1} {_2F_1} \biggl[
\frac{k_i+1}{2};\frac{n-k_0}{2};\frac{k_i+3}{2}; \frac
{(1-Q_{i0}^{-1})(k_i+k_0)}{(1+n)}
\biggr],
\end{eqnarray}
where $_2F_1$ is the standard hypergeometric function [see \citet{AS64}],
and $Q_{i0}=\operatorname{SSE}_i/\operatorname{SSE}_0$ is the ratio of the sum of
squared errors of models $M_i$ and $M_0$.

To implement Bayesian model selection through (\ref{postprob}), one
also needs the prior odds ratios $P_{{j}{0}}$. A recommended objective
Bayesian choice of these odds ratios for the variable selection problem
is $P_{{j}{0}}=k_j!(p-k_j)!/p!$. For extensive discussion and earlier
references, see \citet{ScottBerger09}.

%%-------------------------------------- Acknowledgements
%
%This paper was supported in part by the Spanish Ministry of Education
%and Science under grant MTM2010-19528, and by USA National Science
%Foundation Grants DMS-0635449, DMS-0757549-001, and DMS-1007773.

%--------------------------------------Bibliography
%

%-------------------------------------Appendix

\begin{appendix}

\section*{Appendix}\label{proofs}

\subsection{Approximations to improper priors}\label{app1}

%le1 #&#
\begin{lem}\label{approx-improper}
Consider $\pi_i(\bolds{\alpha}) = c_i \psi_i(\bolds
{\alpha})$, where $\psi_i (\bolds{\alpha})$ increases
monotonically in $i$
to $\pi(\bolds{\alpha})$ and $c_i = 1/\int\psi_i(\bolds
{\alpha}) \,d\bolds{\alpha} <
\infty$. Then,
if $\int f_l(\mathbf{y} \mid\bolds{\alpha})   \pi
(\bolds{\alpha})  \,d\bolds{\alpha} < \infty
$ for all densities $ f_l(\mathbf{y} \mid\bolds{\alpha})$,
\[
\lim_{i \rightarrow\infty} \frac{\int f_l(\mathbf{y} \mid
\bolds{\alpha})  \pi_i(\bolds{\alpha})   \,d\bolds{\alpha}}{\int f_{l^{\prime
}}(\mathbf{y} \mid\bolds{\alpha})
\pi_i(\bolds{\alpha})   \,d\bolds{\alpha}} = \frac{\int f_l(\mathbf{y} \mid\bolds{\alpha})  \pi
(\bolds{\alpha})   \,d\bolds{\alpha}}{\int
f_{l^{\prime}}(\mathbf{y} \mid\bolds{\alpha})
\pi(\bolds{\alpha})   \,d\bolds{\alpha}} .
\]
\end{lem}
\begin{pf}
\[
\frac{\int f_l(\mathbf{y} \mid\bolds{\alpha})  \pi_i(\bolds{\alpha})   \,d\bolds{\alpha}
}{\int f_{l^{\prime}}(\mathbf{y} \mid\bolds{\alpha})
\pi_i(\bolds{\alpha})   \,d\bolds{\alpha}} = \frac
{\int f_l(\mathbf{y} \mid\bolds{\alpha})  \psi_i(\bolds{\alpha})
\,d\bolds{\alpha}}{\int f_{l^{\prime}}(\mathbf{y} \mid
\bolds{\alpha})  \psi_i(\bolds{\alpha}
)   \,d\bolds{\alpha}} \longrightarrow \frac{\int f_l(\mathbf{y} \mid\bolds{\alpha})  \pi
(\bolds{\alpha})   \,d\bolds{\alpha}}{\int
f_{l^{\prime}}
(\mathbf{y} \mid\bolds{\alpha})  \pi(\bolds{\alpha
})   \,d\bolds{\alpha}}
\]
by the monotone convergence theorem.
\end{pf}
Thus common proper priors can be used to approximate common improper
priors and, as the approximation improves, the Bayes factors for the
proper priors converge to the Bayes factor for the improper prior; this
is why Bayesians have always said that it is not illogical to use an
improper prior for a common parameter~$\bolds{\alpha}$ in
computing a Bayes factor.
It is interesting that no conditions are needed in the lemma, except
that the marginal likelihoods exist for the improper prior, which is
clearly needed for the Bayes factor to even be defined for the improper prior.

%-------------------------------------------------------
\subsection{\texorpdfstring{Proof of Proposition \protect\ref{piRSt}}{Proof of Proposition 1}}\label{app2}
This proof requires the following lemma:
%
%le2 #&#
\begin{lem}\label{lempiRSt}
If $m>1$, $p>0$, $a>0$ and $k \geq1$, then
\[
\lim_{z\rightarrow\infty} z^{a+k} \int_0^1
\lambda^{a-1} \biggl(\frac
{\lambda}{m-\lambda} \biggr)^k
e^{-({\lambda}/{(m-\lambda)})\cdot
p\cdot
z}\,d\lambda=m^a \Gamma(a+k) p^{-(a+k)}.
\]
\end{lem}

\begin{pf}
For $0<\varepsilon<1$, write
\begin{eqnarray}
&&\lim_{z\rightarrow\infty}\int_0^1
z^{a+k} \lambda^{a-1} \biggl(\frac{\lambda}{m-\lambda}
\biggr)^k e^{-({\lambda}/{(m-\lambda)})\cdot p\cdot z} \,d\lambda
\nonumber\\
&&\qquad=\lim_{z\rightarrow\infty}\int_0^\varepsilon
z^{a+k}\lambda^{a-1} \biggl(\frac{\lambda}{m-\lambda}
\biggr)^k e^{-({\lambda}/{(m-\lambda)})\cdot p\cdot z} \,d\lambda
\\
&&\qquad\quad{}+\lim_{z\rightarrow\infty}\int_\varepsilon^1
z^{a+k}\lambda^{a-1} \biggl(\frac{\lambda}{m-\lambda}
\biggr)^k e^{-({\lambda}/{(m-\lambda)})\cdot p\cdot z} \,d\lambda
.\nonumber
\end{eqnarray}
Note that
\[
\lim_{z\rightarrow\infty} z^{a+k} \lambda^{a-1} \biggl(
\frac
{\lambda
}{m-\lambda} \biggr)^k e^{-({\lambda}/{(m-\lambda)})\cdot p\cdot z}=0,
\]
and the integrand in the last integral in \eqref{eqAppPr1Lem} is
uniformly bounded over $\lambda$ and $z$. It follows from the dominated
convergence theorem that the last term is zero, so that
%
%e28 #&#
\begin{eqnarray}
\label{eqAppPr1Lem}
&&\lim_{z\rightarrow\infty}\int_0^1
z^{a+k}\lambda^{a-1} \biggl(\frac{\lambda}{m-\lambda}
\biggr)^k e^{-({\lambda}/{(m-\lambda)})\cdot p\cdot z} \,d\lambda
\nonumber
\\[-8pt]
\\[-8pt]
\nonumber
&&\qquad=\lim_{z\rightarrow\infty}\int_0^\varepsilon
z^{a+k}\lambda^{a-1} \biggl(\frac{\lambda}{m-\lambda}
\biggr)^k e^{-({\lambda}/{(m-\lambda)})\cdot p\cdot z} \,d\lambda.
\end{eqnarray}
Next, make the change of variables $t=\lambda/(m-\lambda)$ to get
\[
\int_0^\varepsilon \lambda^{a-1} \biggl(
\frac{\lambda}{m-\lambda} \biggr)^k e^{-({\lambda}/{(m-\lambda)})\cdot p\cdot z} \,d\lambda=m^a
\int_0^{{\varepsilon}/ {(m-\varepsilon)}}
\frac{t^{k+a-1}}{(1+t)^{a+1}} e^{-t\cdot p\cdot z} \,dt.
\]

To bound the integral of interest notice that, for $t\in(0,\varepsilon
/(m-\varepsilon))$,
%
%e29 #&#
\begin{equation}
\label{eqineqaux} \frac{1}{(1+\varepsilon/(m-\varepsilon))^{a+1}}\le\frac
{1}{(1+t)^{a+1}}\le1.
\end{equation}
By integrating $t$ out from \eqref{eqineqaux} and multiplying the
result by $z^{(a+k)}$, we get both an upper and a lower bound for the
integral of interest, namely
%
%e30 #&#
\begin{eqnarray}
\label{eqineq}
&&\frac{m^a  p^{-(a+k)}  (\Gamma(a+k)-\Gamma(a+k,
({\varepsilon}/{(m-\varepsilon)})  p z) )}{(1+\varepsilon/(m-\varepsilon
))^{a+1}}
\nonumber\\[-2pt]
&&\qquad\leq \lim_{z\rightarrow\infty} m^a \int_0^{{\varepsilon}/{(m-\varepsilon)}} z^{a+k} \frac
{t^{k+a-1}}{(1+t)^{a+1}} e^{-t\cdot p\cdot z}
\,dt
\\[-2pt]
&&\qquad\leq m^a p^{-(a+k)} \biggl(\Gamma(a+k)-\Gamma
\biggl(a+k, \frac{\varepsilon
}{m-\varepsilon} p z\biggr) \biggr) ,\nonumber
\end{eqnarray}
where $\Gamma(\nu_1, \nu_2)$ is the incomplete gamma function,
\[
\Gamma(\nu_1, \nu_2)=\int_{\nu_2}^\infty
t^{\nu_1-1}e^{-t} \,dt,
\]
which goes to zero as $\nu_2$ goes to infinity.

Taking limits in~\ref{eqineq} as $z\rightarrow\infty$ gives
\begin{eqnarray*}
\frac{m^a  p^{-(a+k)} \Gamma(a+k)}{(1+\varepsilon/(m-\varepsilon
))^{a+1}} &\leq& \lim_{z\rightarrow\infty} m^a \int
_0^{{\varepsilon} /{(m-\varepsilon)}}
z^{a+k} \frac
{t^{k+a-1}}{(1+t)^{a+1}} e^{-t\cdot p\cdot z} \,dt
\\[-2pt]
&\leq& m^a p^{-(a+k)} \Gamma(a+k).
\end{eqnarray*}
The result follows from \eqref{eqAppPr1Lem} and the fact that the
upper and lower bound are equal as $\varepsilon$ goes to 0.
% , the result follows
% \[
% \lim_{z\rightarrow\infty} m^a z^{a+k} \int_0^\frac{\varepsilon}{m-
%{(1+t)^{a+1}} e^{-t p z} \,dt= m^a  p^{-(a+k)} \Gamma(a+k),
% \]
%proving the result.
\end{pf}
Continuing with the proof of Proposition~\ref{piRSt}, we remove the
subindex $i$ for simplicity in notation.
Since the multivariate Student density can be written as
\begin{eqnarray*}
 \St_k\bigl(\bolds{\beta}\mid\mathbf{0}, \mathbf {C}^{*},
2a\bigr)&=&\frac{\Gamma(a+k/2)}{\Gamma
(a)}(2\pi)^{-k/2} \bigl(\bigl(a\Gamma(a)
\bigr)^{1/a}\sigma^2\rho (b+n) \bigr)^{-k/2}
\\[-2pt]
&&{}\times \bigl|\mathbf{V}^t\mathbf{V}\bigr|^{1/2} \bigl(1+ \bigl(2\bigl(a
\Gamma(a)\bigr)^{1/a}\rho \sigma^2(b+n) \bigr)^{-1}
\|\beta\|^2 \bigr)^{-(a+k/2)},
\end{eqnarray*}
it can be easily shown that
\begin{eqnarray*}
&&\lim_{\|\beta\|^2\rightarrow\infty}\frac{\St_k(\bolds
{\beta}\mid\mathbf{0},  \mathbf{C}^{*}, 2a)} {
\Gamma(a+k/2)(2\pi)^{-k/2} a  (\sigma^2\rho (b+n) )^a
|\mathbf{V}^t\mathbf{V}|^{1/2} 2^{a+k/2}  (\|
\bolds{\beta}\|^2 )^{-(a+k/2)}}
\\[-2pt]
&&\qquad= \biggl(2\bigl(a\Gamma(a)\bigr)^{1/a}\rho\sigma^2(b+n)\\[-2pt]
&&\qquad\hspace*{14pt}{}\times
\lim_{\|\beta\|
^2\rightarrow
\infty} \frac{1+ (2(a\Gamma(a))^{1/a}\rho\sigma^2(b+n)
)^{-1}\|
\beta\|^2}{\|\beta\|^2} \biggr)^{-(a+k/2)}=1.
\end{eqnarray*}
It then follows that
\begin{eqnarray*}
&&\lim_{\|\beta\|^2\rightarrow\infty} \frac{\pi^R(\bolds
{\beta}\mid\bolds{\beta}_{0},  \sigma)}{\St_k(\bolds{\beta}\mid\mathbf{0},
\mathbf{C}^{*}, 2a)}\\
&&\qquad=\frac{(2\sigma^2)^{-(a+k/2)} b^{-k/2}}{\Gamma(a+k/2)  (\rho(b+n) )^a}
\\
&&\qquad\quad{}\times\lim_{\|\beta\|^2\rightarrow\infty} \bigl(\|\bolds{\beta }\|^2
\bigr)^{a+k/2} \int_0^1
\lambda^{a-1} \biggl(\frac{\lambda}{m-\lambda} \biggr)^{k/2}
e^{-{\lambda
}/{(m-\lambda)
} p \|\beta\|^2} \,d\lambda,
\end{eqnarray*}
where $m=(\rho  (b+ n))/b$ and $p=1/(2\sigma^2 b)$. Since $\rho
>b/(b+n)$ and $m>1$, we can apply Lemma~\ref{lempiRSt}, and the
result follows.

%------------------------------------------------------
\subsection{\texorpdfstring{Proof of Result \protect\ref{resIcrit}}{Proof of Result 1}}\label{app3}

To apply invariance, let $\bolds{\theta}= (\bolds{\beta
}_0, \sigma, \bolds{\beta}_i,
i)$ denote the parameter indexing all the models,
and consider the location-scale group defined by $g=(c, \mathbf
{b}) \in
G_0=(0,\infty) \times\mathcal{ R}^{k_0}$ acting on
$\mathbf{y}$ through the transformation $\tilde{\mathbf
{y}}=c  \mathbf{y}+\mathbf{X}_0\mathbf{b}$.
It can be easily seen that
$\tilde{\mathbf{y}}\sim{f}(\cdot\mid\bolds{\theta}^*)$,
where $\bolds{\theta}^{*}=(\bolds{\beta}_{0}^{*}, \sigma^{*}, \bolds{\beta}_i^*,i^*)$
with $\bolds{\beta}_{0}^*=\mathbf{b}+c\bolds{\beta
}_{0}$, $\sigma^*=c\sigma$, $\bolds{\beta}_i^* = c \bolds{\beta}_i$ and $i^*=i$, so that the
transformed model has exactly the same structure as the original model.
The Invariance-criterion thus says that
the prior $\pi_i(\bolds{\beta}_i \mid\bolds{\beta}_0,
\sigma)$ must be such that
the marginal models in (\ref{eqinv-marginals})
are invariant with respect to the group action, so that (keeping to the
notation above)
\[
f\bigl( \tilde{\mathbf{y}} \mid\bolds{\beta}_0^*, \sigma^*,i^*
\bigr) = \int\N_n\bigl(\tilde{\mathbf{y}} \mid\mathbf
{X}_0\bolds{\beta}_{0}^* +\mathbf{X}\bolds{\beta }^*,
\bigl(\sigma^*\bigr)^2\mathbf{I}\bigr) \pi_i\bigl(\bolds{
\beta}^* \mid \bolds{\beta}_{0}^*,\sigma^*\bigr) \,d\bolds{\beta}^* ,
\]
the fact that $\pi_i(\cdot\mid\cdot,\cdot)$ must have the same
functional form as in the original parameterization following from the
completeness of $\N_n(\tilde{\mathbf{y}} \mid\mathbf
{X}_0\bolds{\beta}_{0}^* +\mathbf{X}\bolds{\beta}^*,(\sigma^*)^2\mathbf{I})$, given that the design matrix is of
full rank.
But one can also compute
$f( \tilde{\mathbf{y}} \mid\bolds{\beta}_0^*, \sigma^*,i^*)$ by change of
variables from the original density, yielding
\begin{eqnarray*}
&&f\bigl( \tilde{\mathbf{y}} \mid\bolds{\beta}_0^*, \sigma^*,i^*
\bigr)
\\
&&\qquad=\int\N_n\bigl(\tilde{\mathbf{y}} \mid\mathbf{X}_0
\bolds {\beta}_{0}^* +\mathbf{X}\bolds{\beta}^*,\bigl(\sigma^*
\bigr)^2\mathbf{I}\bigr) \pi_i\bigl(\bolds{\beta}^*/c\mid
\bigl(\bolds{\beta}_{0}^*-\mathbf{b}\bigr)/c,\sigma^*/c\bigr)
c^{-k_0} \,d\bolds{\beta}^* .
\end{eqnarray*}
Again using the completeness of the normal density, these two
expressions can be equal only if
\[
\pi_i\bigl(\bolds{\beta}^* \mid\bolds{\beta}_{0}^*,
\sigma^*\bigr) = \pi_i\bigl(\bolds{\beta}^*/c\mid\bigl(\bolds{\beta
}_{0}^*-\mathbf{b}\bigr)/c,\sigma^*/c\bigr) c^{-k_0} .
\]
This condition is satisfied by the conditional prior in
(\ref{eqscale}).\vadjust{\goodbreak}

With respect to the only ``if'' part of the proof, note that for the
particular transformation in $G_0$ given by $\mathbf
{b}=\bolds{\beta}_{0}$ and
$c=\sigma^*$, the above condition becomes
\[
\pi_i\bigl(\bolds{\beta}^* \mid\bolds{\beta}_{0}^*,
\sigma^*\bigr) =\sigma^{-k_0} \pi\bigl(\bolds{\beta}^*/\sigma^* \mid
\mathbf{0}, 1\bigr),
\]
proving that being of the form in (\ref{eqscale}) is also a necessary
condition.

%-----------------------------------------------
\subsection{\texorpdfstring{Proof of Result \protect\ref{resPM2}}{Proof of Result 2}}\label{app4}

With the use of the full conditional for $\bolds{\beta}_i$
associated with
this prior, the integrated models can be alternatively expressed as
\[
M_i^I\dvtx \mathbf{Y}^*=\mathbf{X}_0
\bolds{\beta }_0+\sigma\bolds{\varepsilon},
\]
where $\bolds{\varepsilon}\sim f_i^I(\mathbf{u})$, given by
\[
f_i^I(\mathbf{u})=\int\N_n(\mathbf{u}\mid
\mathbf {X}_i\mathbf{t},\mathbf{I}) h_i(\mathbf{t}) \,d
\mathbf{t}\qquad (i>0) \quad\mbox{and}\quad f_0^I(\mathbf{u})=
\N_n(\mathbf {u}\mid\mathbf{0} ,\mathbf{I}).
\]
This model selection problem was explicitly studied in \citet
{BerPerVar98}, where it was shown that the minimal sample size
associated with the right-Haar prior for $(\bolds{\beta
}_0,\sigma)$ is
$n^*_i=k_0+1$, and that it is sufficient for exact predictive matching
for $f_i^I(\cdot)$ [or, equivalently, $h_i(\cdot)$] to be symmetric
about the origin.

%--------------------------------------------------
\subsection{\texorpdfstring{Proof of Result \protect\ref{resPM1}}{Proof of Result 3}}\label{app5}

It is convenient to work in terms of orthogonal parameters, so, for
each model $M_{i}$, define $\bolds{\gamma}=\bolds{\beta
}_0+(\mathbf{X}_0^t\mathbf{X}_0)^{-1}\mathbf
{X}_0^t\mathbf{X}_i \bolds{\beta}_i$; this will be
``common'' to all
models and orthogonal to $\bolds{\beta}_i$ in each model $M_i$,
which can be
written in the new parameterization as $\mathbf{y} \sim\N_{n}(\mathbf{y}\mid\mathbf{X}_{0}\bolds{\gamma
}+\mathbf{V}_{i}\bolds{\beta}_{i},  \sigma^{2}\mathbf
{I}_{n})$.
Consider a scale mixture of normals prior of the form
\[
\pi(\bolds{\beta}_{i}\mid\bolds{\gamma}, \sigma)=\pi (\bolds{
\beta}_{i}\mid\sigma)=\int_{0}^{\infty}
\N_{k_{i}}\bigl(\bolds{\beta}_{i}\mid\mathbf {0}, g
\sigma^{2}\mathbf{A}_i\bigr) h(g) \,dg .
\]
Noting that the right-Haar prior for $(\bolds{\alpha}, \sigma)$
transforms
into the same prior ($1/\sigma$) for $(\bolds{\gamma}, \sigma
)$, it follows
that the marginal likelihood under model $M_{i}$ is
\begin{eqnarray*}
m_{i}(\mathbf{y}) &=&\int\N_{n}\bigl(\mathbf{y}\mid\mathbf
{X}_{0}\bolds{\gamma}+\mathbf{V}_{i}\bolds{
\beta}_{i}, \sigma^{2}\mathbf{I}_{n}\bigr)
\sigma^{-1} \pi_{i}(\bolds{\beta}_{i}\mid\bolds{
\gamma}, \sigma) \,d(\bolds{\beta}_{i}, \bolds{\gamma}, \sigma)
\\
&=& \int_{0}^{\infty} \int\N_{n}\bigl(
\mathbf{y}\mid\mathbf {X}_{0}\bolds{\gamma}+\mathbf{V}_{i}
\bolds{\beta }_{i}, \sigma^{2}\mathbf{I}_{n}
\bigr) \sigma^{-1} \\
&&\hspace*{6pt}\qquad{}\times \N_{k_{i}}\bigl(\bolds{\beta}_{i}
\mid\mathbf{0}, g \sigma^{2}\mathbf{A}_i\bigr)
 h(g) \,d(\bolds{\beta}_{i}, \bolds{\gamma}, \sigma) \,dg .
\end{eqnarray*}
Using the fact that $\mathbf{y}^{t}\mathbf{V}_{i}(\mathbf
{V}_{i}^{t}\mathbf{V}_{i})^{-1}\mathbf{V}_{i}^{t}\mathbf
{y} = \operatorname{SSE}_{0}$ for any sample of size $n=k_{i}+k_{0}$ and
integrating out $\bolds{\gamma}$, $\bolds{\beta}_{i}$ and
$\sigma$ yields
\begin{eqnarray*}
m_{i}(\mathbf{y})&=&\int_{0}^{\infty
}\bigl|
\mathbf{X}_{0}^{t}\mathbf{X}_{0}\bigr|^{-1/2}
\frac{\pi^{-k_{i}/2}|(\mathbf{V}_{i}^{t}\mathbf
{V}_{i})^{-1}|^{1/2}}{2| (\mathbf{V}_{i}^{t}\mathbf
{V}_{i})^{-1}+g\mathbf{A}_i|^{1/2}}
\\
&&\hspace*{18pt}{}\times \bigl(\hat{\bolds{\beta}}{}^{t}_{i}\bigl[\bigl(\mathbf
{V}_{i}^{t}\mathbf{V}_{i}\bigr)^{-1}+g
\mathbf{A}_i\bigr]^{-1}\hat {\bolds{\beta}}_{i}
\bigr)^{-k_{i}/2}\Gamma \biggl(\frac
{k_{i}}{2} \biggr)h(g) \,d(g).
\end{eqnarray*}
For the robust prior, $\mathbf{A}_i =(\mathbf
{V}_{i}^{t}\mathbf{V}_{i})^{-1}$, and it
follows that
\[
m_{i}(\mathbf{y})=\frac{1}{2}\bigl|\mathbf{X}_{0}^{t}
\mathbf {X}_{0}\bigr|^{-1/2}\pi^{-k_{i}/2} \Gamma \biggl(
\frac{k_{i}}{2} \biggr) (\operatorname{SSE}_{0})^{-k_i/2} ,
\]
which is the same for all models of dimension $k_i$, establishing that
the robust prior is dimension predictive matching for sample sizes
$k_0+k_i$. Furthermore, this last expression equals $m_{0}(\mathbf
{y})$ (see
Appendix~\ref{app6}), establishing that the robust prior is null predictive
matching for samples of size $k_0+k_i$. [Note that this result would
hold for any proper choice of $h(g)$, not just that for the robust
prior.]

To see that null predictive matching does not occur if $\mathbf
{A}_i$ is not
a multiple of $(\mathbf{V}_{i}^{t}\mathbf{V}_{i})^{-1}$, note
that the expression
to be established for null predictive matching is (eliminating
multiplicative constants)
\begin{eqnarray*}
0&=&\int_{0}^{\infty} \biggl(\frac{|(\mathbf
{V}_{i}^{t}\mathbf{V}_{i})^{-1}|^{1/2}(\hat{\bolds{\beta}}{}^{t}_{i}[(\mathbf{V}_{i}^{t}\mathbf{V}_{i})^{-1}+g\mathbf
{A}_i]^{-1}\hat{\bolds{\beta}}_{i})^{-k_{i}/2}}{| (\mathbf
{V}_{i}^{t}\mathbf{V}_{i})^{-1}+g\mathbf{A}_i|^{1/2}}\\
&&\hspace*{124pt}\qquad{} -\bigl(
\hat{\bolds{\beta}}{}^t_{i} \mathbf{V}_{i}^{t}
\mathbf {V}_{i}\hat{\bolds{\beta}}_{i}\bigr)^{-k_{i}/2}
\biggr)
h(g) \,d(g).
\end{eqnarray*}
Since $(\mathbf{V}_{i}^{t}\mathbf{V}_{i})^{-1}$ and
$\mathbf{A}_i$ are positive definite,
there is a matrix $\mathbf{B}$ such that $\mathbf
{B}^t(\mathbf{V}_{i}^{t}\mathbf{V}_{i})^{-1}\mathbf{B} =
\mathbf{I}$ and $\mathbf{B}^t \mathbf{A}_i \mathbf
{B} = \mathbf{D}$, with $ \mathbf{D}$ being a
diagonal matrix with diagonal elements $d_i$. Also defining
$\mathbf{W} = {\mathbf{B}}^{t} \hat{\bolds{\beta
}}_{i}$, it follows that the above
expression can be written
\[
0 = \int_{0}^{\infty} \biggl(\frac{(W^t[\mathbf{I}+g\mathbf
{D}]^{-1}W)^{-k_{i}/2}}{|\mathbf{I}+g\mathbf{D}|^{1/2}} -
\bigl(|W|^2\bigr)^{-k_{i}/2} \biggr) h(g) \,d(g) .
\]
Let $d_j$ be the largest diagonal element, and choose $W$ to be the
unit vector in coordinate $j$. Then the above expression becomes
\[
0 = \int_{0}^{\infty} \biggl(\frac{(1+g d_j)^{k_{i}/2}}{ \prod_{l=1}^{k_i}(1+gd_i)^{1/2}} -1
\biggr) h(g) \,d(g) .
\]
But the integrand is clearly greater than 0, unless all $d_i$ are equal
which is equivalent to the statement that $\mathbf{A}_i$ is a
multiple of
$(\mathbf{V}_{i}^{t}\mathbf{V}_{i})^{-1}$.

%----------------------------------------------------------
\subsection{\texorpdfstring{Computation of the Bayes factor in (\protect\ref{eqBFR})}
{Computation of the Bayes factor in (16)}}\label{app6}

%pr2 #&#
\begin{prop} For any $(a,  b,  \rho_{i})$ satisfying (\ref
{eqg-conditions}) and $n\geq k_{i}+k_{0}$, the prior predictive
distribution for $\mathbf{y}$ under $M_{i}$ using the robust prior is
\[
m^{R}_i(\mathbf{y})=m^{R}_{0}(
\mathbf{y}) Q_{i0}^{-
{(n-k_0)}/{2}} \frac
{2a}{k_i+2a} \bigl[
\rho_{i} (n +b)\bigr]^{-{k_i}/{2}} \operatorname{AP}_{i0},
\]
where
\[
m^{R}_{0}(\mathbf{y})=\frac{1}{2}
\pi^{-{(n-k_{0})}/{2}} \bigl|\mathbf{X}_0^t \mathbf{X}_0\bigr|^{-{1}/{2}}
\Gamma \biggl[\frac{n-k_{0}}{2} \biggr] \operatorname {\operatorname{SSE}}_0^{
 -{(n-k_{0})}/{2}}\vadjust{\goodbreak}
\]
and $\operatorname{AP}_{i}$ defined in (\ref{eqapell}).
Hence the Bayes factor obtained with prior $\pi_i^R$ in \eqref
{eqourprior} can be compactly expressed as in (\ref{eqBFR}).
\end{prop}
\begin{pf}
It is convenient to carry out the proof in the orthogonal
transformation of the parameters as in Appendix~\ref{app5}.
Using standard normal computations, the prior predictive distribution
under $M_0$ is
\begin{eqnarray*}
m^{R}_0(\mathbf{y}) &=&\int_{\mathbb{R}^{k_{0}}}
\int_0^\infty \N_n\bigl(\mathbf{y}\mid
\mathbf{X}_0\bolds{\gamma}, \sigma^2
\mathbf{I}_n\bigr) \frac{1}{\sigma} \,d\bolds{\gamma} \,d\sigma
\\
&=& \frac{1}{2}\pi^{-{(n-k_0)}/{2}} \bigl|\mathbf{X}_0^t
\mathbf{X}_0\bigr|^{-
{1}/{2}} \Gamma \biggl[\frac{n-k_0}{2}
\biggr] \operatorname{SSE}_0^{  -
{(n-k_0)}/{2}} .
\end{eqnarray*}
Integrating out $\bolds{\beta}_i$, $\bolds{\gamma}$ and
$\sigma$, the prior
predictive distribution under $M_i$ is
\begin{eqnarray*}
m^{R}_i(\mathbf{y}) &=&\int\N_n\bigl(
\mathbf{y}\mid\mathbf {X}_0\bolds{\gamma}+\mathbf{V}_i
\bolds{\beta}_i, \sigma^2 \mathbf{I}_n\bigr)
\N_{k_i}\bigl(\bolds{\beta}_i \mid0, \mathbf{B}(\lambda)
\bigr)
 a \lambda^{a-1} \sigma^{-1} \,d(
\bolds{\gamma}, \bolds{\beta}_i, \sigma, \lambda)\\
&=&\frac{1}{2}
\pi^{-{(n-k_0)}/{2}} \bigl|\mathbf {X}_0^t\mathbf{X}_0\bigr| ^{-{1}/{2}}
\Gamma \biggl[\frac
{n-k_0}{2} \biggr]
\\
&&{}\times\int_0^1 a\lambda^{a+{(k_i}/{2})-1}\bigl(
\rho_i(b+n)-(b-1)\lambda \bigr)^{{(n-k_i-k_{0})}/{2}}
\\
&&\hspace*{27pt}{}\times\bigl(\operatorname{SSE}_i \bigl(\rho_i (b+n)-b
\lambda\bigr)+\lambda \operatorname {SSE}_0 \bigr)^{-{(n-k_0)}/{2}} \,d
\lambda ,
\end{eqnarray*}
with $\mathbf{B}(\lambda)=(\lambda^{-1}\rho_{i}(b+n)-b)\sigma^{2}(\mathbf{V}_{i}^{t}\mathbf{V}_{i})^{-1}$. This expression
can be rewritten as
\begin{eqnarray*}
m_{i}^R(\mathbf{y}) &=& a Q_{i0}^{-{(n-k_0)}/{2}}
\bigl(\rho_i (n +b)\bigr)^{-k_i/2} m_{0}^R(
\mathbf{y})
\\
&&{}\times\int_0^1\lambda^{a+({k_i}/{2})-1}
\biggl(1-\frac{b-1}{\rho_i(b+n)}\lambda \biggr)^{{(n-k_i-k_0)}/{2}}\\
&&\hspace*{4pt}\qquad{}\times \biggl(1-
\frac{b-Q_{i0}^{-1}}{\rho_i(b+n)}\lambda \biggr)^{-
{(n-k_0)}/{2}}\,d\lambda ,
\end{eqnarray*}
and the result follows by noting that
\begin{eqnarray*}
\operatorname{AP}_{i}&=&\frac{2a+k_{i}}{2}\int_0^1
\lambda^{a+({k_i}/{2})-1} \biggl(1-\frac{b-1}{\rho_i(b+n)}\lambda
\biggr)^{{(n-k_i-k_0)}/{2}}\\
&&\hspace*{51pt}{}\times
\biggl(1-\frac{b-Q_{i0}^{-1}}{\rho_i(b+n)}\lambda \biggr)^{-
{(n-k_0)}/{2}}\,d\lambda .
\end{eqnarray*}
\upqed\end{pf}

%----------------------------------------------------------
\subsection{\texorpdfstring{Proof of Corollary \protect\ref{resCM2}}{Proof of Corollary 1}}\label{app7}
For the prior in \eqref{eqourprior},
\[
\int_0^\infty(1+g)^{-{k_i}/{2}}
p_{i}^{R}(g) \,dg = \int_{\rho
_i(b+n)-b}^\infty(1+g)^{-{k_i}/{2}}
\frac{a [\rho_{i}(b+n)]^a}{(g+b)^{(a+1)}} \,dg.
\]
The change of variables $z=g-[\rho_i(b+n)-b]$ results in
\begin{eqnarray*}
&&\int_0^\infty(1+g)^{-{k_i}/{2}}
p_{i}^{R}(g) \,dg \\
&&\qquad=\int_0^\infty
\frac{a[\rho_i(b+n)]^a}{[z+\rho_i(b+n)]^{(a+1)}
[1+z+\rho_i(b+n)-b]^{{k_i}/{2}}} \,dz.
\end{eqnarray*}
It is now easy to see that if $\rho_i(b+n)$ goes to $\infty$ with $n$,
this integral vanishes as $n\rightarrow\infty$, satisfying the
condition of Result~\ref{resCM1}.

%-----------------------------------------------------
\subsection{\texorpdfstring{Proof of Result \protect\ref{resCI}}{Proof of Result 7}}\label{app8}
For simplicity, the explicit dependence of $Q_{i0}$ on $\mathbf
{y}_{m}$ will
not be shown in this proof, and $\lim_{m\rightarrow\infty
}Q_{i0}(\mathbf{y}_{m})=0$ will be denoted by $Q_{i0}\rightarrow0$.
The robust Bayes factor can be written as
\begin{eqnarray*}
\label{eqAppPr3-1}
B^R_{i0}&=&a \bigl(
\rho_i(n+b)\bigr)^{-{k_i}/{2}}(Q_{i0})^{-
{(n-k_0)}/{2}}\\
&&{}\times \int
_{0}^{1} \lambda^{a+({k_i}/{2})-1} \biggl[1-
\frac{b-1}{\rho_i(b+n)}\lambda \biggr]^{
{(n-k_i-k_0)}/{2}}
 \\
 &&\hspace*{27pt}{}\times \biggl[1-\frac{b-Q_{i0}^{-1}}{\rho_i(b+n)}\lambda \biggr]^{-{(n-k_0)}/{2}}\,d\lambda
\\
&=& a \bigl(\rho_i(n+b)\bigr)^{-{k_i}/{2}}\int
_{0}^{1} \lambda^{a+({k_i}/{2})-1} \biggl[1-
\frac{b-1}{\rho_i(b+n)}\lambda \biggr]^{{(n-k_i-k_0)}/{2}}\\
&&\hspace*{94pt}{}\times\biggl[Q_{i0} \biggl(1-\frac{b\lambda}{\rho_i(b+n)} \biggr)+\frac{\lambda}{\rho_i(b+n)}
\biggr]^{-
{(n-k_0)}/{2}} \,d\lambda.
\end{eqnarray*}
Note that, since $b >0$, $\rho_i \geq b/(b+n)$ and $0< \lambda< 1$,
\[
\min\biggl\{1,\frac{1}{b}\biggr\} \leq \biggl[1-\frac{b-1}{\rho_i(b+n)}\lambda
\biggr] \leq\max\biggl\{1,\frac{1}{b}\biggr\}
\]
and
\[
\biggl[\frac{\lambda}{\rho_i(b+n)} \biggr] \leq \biggl[Q_{i0} \biggl(1-
\frac
{b\lambda}{\rho_i(b+n)} \biggr)+\frac{\lambda}{\rho_i(b+n)} \biggr] \leq \biggl[Q_{i0}+
\frac{\lambda}{\rho_i(b+n)} \biggr] .
\]
Applying these bounds, it is immediate that
%
%e31 #&#
\begin{eqnarray}
\label{eqinfcon} &&c_1\int_{0}^{1}
\lambda^{a+({k_i}/{2})-1} [c_2 Q_{i0}+\lambda ]^{-{(n-k_0)}/{2}}
\,d\lambda
\nonumber
\\[-8pt]
\\[-8pt]
\nonumber
&&\qquad\leq B^R_{i0} \leq c_3 \int
_{0}^{1} \lambda^{a+({k_i}/{2})-1} [\lambda
]^{-{(n-k_0)}/{2}} \,d\lambda
\end{eqnarray}
for positive constants $c_1$, $c_2$ and $c_3$.

To prove the ``only if'' part of the proposition, note that the last
integral in (\ref{eqinfcon}) is finite if $n <k_i+k_0+2a$. Hence
$B_{i0}$ is bounded by a constant as $Q_{i0} \rightarrow0$, and
information consistency does not hold.\vadjust{\goodbreak}

To prove the ``if'' part of the proposition, make the change of
variables $\lambda^* = \lambda/Q_{i0}$ in the lower bound in (\ref
{eqinfcon}), resulting in the expression
\[
Q_{i0}^{(2a+k_0+k_i-n)/2} c_1\int_{0}^{Q_{i0}^{-1}}
\bigl(\lambda^*\bigr)^{a+
({k_i}/{2})-1} \bigl[c_2 +\lambda^*
\bigr]^{-{(n-k_0)}/{2}} \,d\lambda^* .
\]
If $n > k_i+k_0+2a$, it is clear that this expression goes to infinity
as $Q_{i0} \rightarrow0$ (since the integral itself cannot go to 0).
If $n = k_i+k_0+2a$, the expression becomes
\[
c_1\int_{0}^{Q_{i0}^{-1}} \biggl(
\frac{ \lambda^*}{c_2 +\lambda^*} \biggr)^{a+({k_i}/{2})} \bigl(\lambda^*\bigr)^{-1} \,d
\lambda^* ,
\]
which clearly goes to infinity as as $Q_{i0} \rightarrow0$, completing
the proof.
\end{appendix}

% imsref loaded by akundreckaite, 2012-06-26 11:25:59

%suskaldyti doi

\printaddresses

\end{document}